\documentclass[reqno]{article}
\usepackage{amssymb, amsmath, amsthm, amsfonts, amscd}
\usepackage[english]{babel}
\usepackage{bbm}

\usepackage{graphicx}
\graphicspath{ {./images/} }
\usepackage{epstopdf}

\usepackage[margin=3.cm]{geometry}
\usepackage{breqn}
\usepackage{mathtools}
\usepackage{color}
\usepackage{hyperref}
\hypersetup{
    colorlinks=true,
    citecolor=red,
    filecolor=black,
    linkcolor=blue,
    urlcolor=black
}

\usepackage{booktabs}

\usepackage{xr}
\chardef\bslash=`\\ 

\newtheorem{theorem}{Theorem}[section]
\newtheorem*{theorem*}{Theorem}
\newtheorem{corollary}[theorem]{Corollary}
\newtheorem{lemma}[theorem]{Lemma}
\newtheorem{prop}[theorem]{Proposition}

\newtheorem*{remark}{Remark}

\newcommand{\N}{\mathbb{N}}

\newcommand{\BB}{\mathcal{B}}

\newcommand{\B}{\mathrm{B}}
\newcommand{\F}{\mathrm{F}}
\newcommand{\FF}{\mathbb{F}}
\newcommand{\FFF}{\mathcal{F}}

\newcommand{\LL}{\mathcal{L}}

\newcommand{\Ham}{\mathcal{H}}

\newcommand{\R}{\mathbb{R}}

\newcommand{\T}{\bar{T}}

\newcommand{\EE}{\mathcal{E}}
\newcommand{\EH}{\mathcal{E}_{0}^{0}}
\newcommand{\EI}{\mathcal{E}_{00}}

\newcommand{\xx}{\dot{x}}
\newcommand{\HH}{H^{1}}
\newcommand{\dx}{\dot{x}}
\newcommand{\h}{^{0}_{0}}
\newcommand{\dd}{\mathrm{d}}

\newcommand{\scIh}{C^{1}_{0}(I)}

\def\a{\alpha }
\def\b{\beta }

\def\l{\lambda }

\def\d{\delta}

\def\w{\omega}
\def\e{\varepsilon}
\def\f{\varphi}
\def\.{\cdot }
\def\ra{\rightarrow}
\def\hra{\hookrightarrow}

\title{
\textsc{\textbf{The weak form of the SDOF and MDOF equation of motion, part I: Theory}}\\
\author{Nikolaos Karaliolios$^1$ \and Dimitrios L. Karabalis$^2$}
\date{%
    $^1$email: nkaraliolios@gmail.com\\%
    $^2$Department of Civil Engineering, University of Patras\\[2ex]%
    \today
}
}

\begin{document}

\maketitle

\begin{abstract}
The weak form of the SDOF and MDOF equations of motion are obtained.
The original initial conditions problem is transformed into a boundary
value problem. The boundary value problem is then solved and transformed
back to the initial conditions one. Subsequently, a general method for
obtaining numerical methods using an arbitrary number of linearly independent
approximating functions is outlined.

This is part one of a series of three papers, in the second of which
a numerical method is obtained, using Bernstein polynomials of arbitrarily high
order. The numerical evidence for the convergence of the method will be
presented in the third part paper.

\end{abstract}  

\tableofcontents


\section{The SDOF problem}

The equation of motion of the single-degree-of-freedom (SDOF) problem is
\begin{equation}
m \ddot{x}  + c \dot{x}  + k x  = f 
\end{equation}
where $m>0 $ is the mass, $c \geq 0 $ the damping coefficient,
$k >0 $ the stiffness, and $f \colon \R \ra \R $ is the excitation
force, which is considered to act for finite time, so that it is $0$ outside the
interval $ [0,\bar{T}]$, where $0<\bar{T}<\infty $. In the notation
of the next paragraph,
this is equivalent to $\mathrm{supp} (f) \subset [0,\bar{T}]$.
In the following, $I$ denotes the open interval $(0,\bar{T})$,
and $\bar{I}$ its closure $[0,\bar{T}]$. Mass normalized to $1$ is the
standing assumption throughout this study of the SDOF problem, while
$c$, $k$ and $f$ are considered fixed. The equation thus reaches
the form
\begin{equation} \label{eq SDOF}
\ddot{x}  + c \dot{x} + k x  = f \text{ in } [0,\bar{T}]
\end{equation}
For later use, the operator defining the lhs of the equation is denoted by
\begin{equation} \label{eqdefF}
\FF (x) : x \ra \ddot{x}  + c \dot{x} + k x
\end{equation}
which is a linear operator $H^{s+2} \ra H^{s}$, for all $s \geq -1$ (see \S \ref{secdef} for
the definition of the Sobolev spaces $H^{s}$). The operator $\FF $ associates a smooth
enough displacement function to the excitation force of the corresponding SDOF problem.

The equation is coupled with initial conditions
\begin{equation} \label{eq in cond}
\begin{cases}
x(0) = x_{0} \in \R \\
\dot{x}(0) = \dot{x}_{0} \in \R
\end{cases}
\end{equation}

In the underdamped case, i.e. $c < 2\sqrt{k}$, which is the case of interest, the
homogeneous case of eq. \eqref{eq SDOF}, i.e. for $f(\. ) \equiv 0$, admits
the well known solution (see, e.g., \cite{clough1975dynamics})
\begin{equation} \label{eq sol hom}
x_{h} (t ) = x_{h} (t ;  x_{0} , \dot{x}_{0} ) =
e^{-\xi \w _{n} t} \left[
x_{0} \cos (\w _{d} t) +
\frac{\dot{x}_{0}+ \xi \w _{n} x_{0}}{\w _{d}} \sin (\w _{d} t) 
\right]
\end{equation}
where the subscript $h$ stands for homogeneous; and
\begin{equation} \label{eq damp}
\xi = \frac{c}{2\sqrt{k}} < 1
\end{equation}
is the damping ratio,
\begin{equation} \label{eq undamp eigenfr}
\w _{n} = \sqrt{k}
\end{equation}
is the natural eigenfrequency of the undamped system, and
\begin{equation} \label{eq damp eigenfr}
\w _{d} = \w _{n} \sqrt{1 - \xi ^{2}}
\end{equation}
is the eigenfrequency of the damped system.

If $f $ is, say, $C^{0}$, the traditional theory of ODEs, see e.g.
\cite{HirschSmaleODE}, provides existence and uniqueness of the solution. Another
way of obtaining the existence and uniqueness for the inhomogeneous problem,
admitting existence and uniqueness for the homogeneous one, is by the Duhamel
convolution formula
\begin{equation} \label{eq Duh}
\begin{array}{r@{}l}
x(t) &= \frac{1}{{\omega _d }}\int_0^t {f(\tau )e^{ - \xi \omega _n (t
- \tau )} \sin [\omega _d (t - \tau )]\dd\tau }
+ x_{h} (t ;  x_{0} , \dot{x}_{0} ) \\
&= \frac{1}{{\omega _d }}\int_0^t {f(\tau )e^{ - c (t
- \tau )/2} \sin [\omega _d (t - \tau )]\dd\tau }
+ x_{h} (t ;  x_{0} , \dot{x}_{0} )
\end{array}
\end{equation}
The formula of eq. \eqref{eq Duh} is meaningful even in lower
regularity, for example $f \in L^{\infty}$ or even $f \in L^{1}$, even
though the differential equation, eq. \eqref{eq SDOF}, is not well defined for
such an $f$ due to regularity issues concerning the definition of the acceleration,
$\ddot{x}$, as a well-defined pointwise function.

\section{From initial to boundary conditions}

The linearity of the problem also provides a natural way of
transforming the initial value problem of eqs. \eqref{eq SDOF} and \eqref{eq in cond} to a boundary-value problem. If the goal is to solve the
equation for times $t \in [0,\bar{T}]$ with $\bar{T}<\infty$, then the
operator
\begin{equation} \label{eq bound prob}
\BB = \BB _{x_{0},f,\bar{T}} \colon \xx _{0} \mapsto x(\bar{T})
\end{equation}
can be defined,
where $x(\. )$ is the solution of the eq. \eqref{eq SDOF} with
initial velocity $\xx_{0}$. The initial position $x_{0}$, the
excitation function $f$ and the time $\bar{T}$ are to be considered as
fixed parameters.

The solution of the homogeneous problem, eq. \eqref{eq sol hom},
which is well defined as a $C^{\infty}$ smooth function of time,
together with the Duhamel representation of the solution show that
for a given set of parameters $x_{0}$, $f$ and $\bar{T}$, the mapping
$\BB : \R \ra \R $ is affine linear, i.e. of the form
$\xx _{0} \mapsto \a \xx _{0} + \beta $ for some real constants $\a$ and
$\beta $, and $\xx _{0} \in \R $. The slope $\a$ is non-zero and, consequently,
the mapping is one-to-one if, and only if, $\sin (\w _{d} \bar{T} ) \neq 0$,
which holds for $\bar{T} $ in an open-dense subset of $\R $.\footnote{An open-dense
subset of $\R$ is a union of open intervals, the complement of which is a
sequence of points. The set where $\sin (\w _{d} \bar{T} ) = 0$ is such a
set, which satisfies the additional property of not having finite
accumulation points. The complement of such a set is exceptional in
any relevant sense.} This is the case where no stationary waves of the eigenfrequency
of the SDOF system fit exactly in the interval $\bar{I}$. In terms of
Control Theory, this is the case where the final displacement can be
controlled by the initial velocity. The slope $\a $ and the
constant $\beta $ depend in a linear way on the parameters $x_{0}$ and
$f $, for fixed $\bar{T} $. Therefore, for a typical $\bar{T} $, the initial
value problem is equivalent to a unique boundary condition problem
\begin{equation} \label{eq bound cond}
\begin{cases}
x(0) = x_{0} \in \R \\
x(\bar{T}) = x_{\bar{T}} = \BB _{x_{0},f,\bar{T}} (\xx _{0}) \in \R
\end{cases}
\end{equation}
As in the underdamped case, a simple calculation can show that for
overdamped systems this operator is always one-to-one, independent of $\bar{T}$.

Using this notation, the goal of this paper is to obtain the weak
formulation of the initial value problem of eq. \eqref{eq SDOF} with
the initial conditions of eq. \eqref{eq in cond}. This amounts to
stating and proving theorem \eqref{thm weak form SDOF}. Subsequently, its
direct generalization to the multi-degree-of-freedom (MDOF) case, theorem
\eqref{thm weak form MDOF}, is stated below and proved.

\section{Definitions} \label{secdef}

The following material is recalled and defined for reasons of
completeness of the article. The reader unfamiliar with the notions
defined below can refer to \cite{RudMathAn} and \cite{RudR&C} for the
mathematical and real analysis related concepts, and to
\cite{BrezisAnFonc} for a general background in functional analysis, the
definition of Sobolev spaces and some applications on the
Sturm-Liouville problem which can serve as an introduction to this
work.

The open interval $I = (0, \bar{T})$ with $\bar{T}< \infty$ is to be
considered fixed, and $\mu $ is the standard Lebesgue measure on $\R$ and its
restriction on $I$.
The space of $C^{s}$ smooth real functions on
$I$ will be denoted by $C^{s}(I,\R )$, or simply by $C^{s}$, for
$0\leq s \leq \infty$.
The space of $C^{s}$ smooth functions on the closed interval
$\bar{I} = [0,\bar{T}]$ will be denoted by $C^{s}(\bar{I},\R )$.
The support of a function $f$, denoted by $\mathrm{supp}f$, is the
closure\footnote{The closure of an open or semi-open interval is the
closed interval with the same endpoints. For a more general
definition of the closure of a set, see \cite{RudMathAn}.}
of the set $\{ x \in I, f(x) \neq 0 \} $. Smooth functions
whose support is contained in $I$ will be denoted by $C^{s}_{0}$.
Informally, such functions, when used as test functions, do not see
(or measure) the behavior of the tested object at the endpoints of the interval
$I$, since the values of the object close enough to the endpoints are
indifferent, as they are integrated against $0$.

The space $L^{2}(I)$ is the Hilbert space of Lebesgue square integrable
functions on $I$ with its standard scalar product
\begin{equation}
\langle f,g \rangle _{L^{2}} = \int _{I}fg \dd \mu
\end{equation}
Another, equivalent, scalar product will be defined in eq. \eqref{eq def in prod c},
which is more adapted to the framework of this article. The space
$L^{2}(I)$ is the Cauchy completion of $C^{\infty}$ for the $L^{2}$
norm
\begin{equation}
\| f \|_{L^{2}} = \left( \int _{I} f^{2} \dd \mu \right) ^{1/2}
\end{equation}
The space of essentially bounded functions on $I$ will be denoted
by $L^{\infty}$, and the space of Lebesgue integrable functions
will be denoted by $L^{1}$.

The Sobolev space $H^{1} = H^{1}(I)$ can be defined as the completion of
$C^{\infty}(\bar{I},\R )$ for the norm induced by the scalar product
\begin{equation}
\langle f,g \rangle _{H^{1}} = \int _{I}fg \dd \mu +
\int _{I}f'g' \dd \mu
\end{equation}
while the Sobolev space $H^{1}_{0} \subset H^{1}$ as the completion of
$C^{\infty}_{0}(I,\R )$ for the same norm. Informally, the functions in
$H^{1}$ have "one derivative in $L^{2}$". An alternative
definition for $H^{1}$ is
\begin{equation}
\{ u \in L^{2} , \exists v \in L^{2} ,  \int _{I} v \f \dd \mu = -\int _{I} u \f ' \dd \mu , \forall
\f \in C^{\infty}_{0}(I,\R) \}
\end{equation}
where $\f $ are test functions. The function $v$ is then by definition
the weak derivative of $u$ and is denoted by $u'$ or by $\dot{u}$ when
the variable is time. The Sobolev
injection theorem states that, in fact, $H^{1} \hra C^{0}(\bar{I})$
and thus functions in $H^{1}_{0}$ can be defined as
\begin{equation}
\{ u \in H^{1} , u(0)=u(\bar{T})=0 \}
\end{equation}
In other words, functions in $H^{1}_{0}$ satisfy homogeneous boundary
conditions.

The Sobolev spaces of higher regularity $H^{s}$ are defined as
the completions of $C^{\infty}(\bar{I},\R )$ for the respective norm
associated to the scalar product
\begin{equation}
\langle f,g \rangle _{H^{s}} = \int _{I}fg \dd \mu +
\int _{I}f^{(s)}g^{(s)} \dd \mu
\end{equation}
and the spaces $H^{s}_{0}$ as the completions of $C^{\infty}_{0}(I,\R )$
for the same respective norms. Functions in $H^{s}$ have $s-1$
continuous derivatives and "their derivative of order $s$ is in
$L^{2}$", while functions in $H^{s}_{0}$ satisfy the homogeneous
boundary conditions
\begin{equation}
u^{(i)}(0) = u^{(i)}(T) = 0, 0 \leq i \leq s-1
\end{equation}

Finally, let a function $f \colon H \ra \R$, where $H $ is, say, a
Hilbert space. The directional, or G\^{a}teaux, derivative of $f$ at
the point $u \in H$ is defined as the mapping of $H$ into itself
given by
\begin{equation}
v \mapsto \lim _{t \ra 0} \frac{f(u+tv) - f(u)}{t\|v\|}
\end{equation}
provided that the limit exists for all $v \in H$, and provided that
the resulting map is a continuous linear selfmapping of $H$.

\section{The weak formulation of the SDOF problem}

The following calculations and transformations of eq. \eqref{eq SDOF}
constitute the preparation for the statement of the weak formulation
of the SDOF problem, i.e. of theorem \ref{thm weak form SDOF}.

Multiplication of eq. \eqref{eq SDOF} by $e^{ct}$ and use of the formula
for the derivative of a product yields
\begin{equation}
\frac{\dd}{\dd t}\left[e^{ct}\xx (t) \right] +  k e^{ct}x(t) = e^{ct}f(t)
\end{equation}
Multiplication by $\f \in C^{\infty}_{0}(I)$ and integration by parts
then gives
\begin{equation} \label{eq. a}
\left[e^{ct}\xx (t)\f (t) \right]_{0}^{\bar{T}} -
\int _{0}^{\bar{T}} e^{ct} \xx (t) \dot{\f} (t) \dd t +
k \int _{0}^{\bar{T}} e^{ct}x(t)\f (t)  \dd t =
\int _{0}^{\bar{T}} e^{ct}f(t)\f (t)  \dd t
\end{equation}
Since $\f (0) = \f (\bar{T}) = 0$, the first term of eq. \eqref{eq. a} vanishes, so that
the equation
\begin{equation} \label{eq weak form interm}
-\int _{0}^{\bar{T}} e^{ct} \xx (t) \dot{\f} (t) \dd t +
k \int _{0}^{\bar{T}} e^{ct}x(t)\f (t)  \dd t =
\int _{0}^{\bar{T}} e^{ct}f(t)\f (t)  \dd t
\end{equation}
is satisfied for every $\f \in C^{\infty}_{0}(I)$. Upon introduction of the space $L^{2}_{c}(I) $, consisting of the square-integrable
functions on $I$ and its inner product
\begin{equation} \label{eq def in prod c}
\langle
\f (\. ) , \psi (\. )
\rangle =
\langle
\f (\. ) , \psi (\. )
\rangle _{c} =
\int _{0}^{\bar{T}} e^{ct}  \f (t)\psi (t) \dd t,
\end{equation}
where $c$ will be suppressed from the notation whenever its value is clear from the context,
eq. \eqref{eq weak form interm} can be written in a more compact form
\begin{equation} \label{eq weak form non hom}
- \langle
\xx (\. ) , \dot{\f} (\. )
\rangle 
+ k 
\langle
x (\. ) , \f (\. )
\rangle =
\langle
f (\. ) , \f (\. )
\rangle , \forall \f \in C^{\infty}_{0}(I)
\end{equation}

Assuming boundary conditions as in eq. \eqref{eq bound cond}, two
$C^{\infty}$ functions $v_{0}, v_{\bar{T}}\colon [0,\bar{T}] \ra \R $ can be
chosen and fixed such that
\begin{equation}
\begin{cases}
v_{0}(0) = v_{\bar{T}}(\bar{T}) = 1 \\
v_{0}(\bar{T}) = v_{\bar{T}}(0) = 0
\end{cases}
\end{equation}
Thus, if $x(\. ) $ solves eq. \eqref{eq SDOF} with boundary
conditions those of eq. \eqref{eq bound cond}, then the function
\begin{equation} \label{eqdef u}
u(\. ) = x(\. ) - x_{0} v_{0}(\. ) - x_{\bar{T}} v_{\bar{T}}(\. )
\end{equation}
satisfies homogeneous boundary conditions $u(0)= u(\bar{T}) =0$.
Calling
\begin{equation}
v_{0,\bar{T}}(\. ) = x_{0} v_{0}(\. ) + x_{\bar{T}} v_{\bar{T}}(\. )
\end{equation}
and
\begin{equation}
u(\. ) = x(\. ) - v_{0,\bar{T}}(\. ) \in H^{1}_{0}
\end{equation}
so that
\begin{equation}
 x(\. ) = u(\. )  + v_{0,\bar{T}}(\. )
\end{equation}
eq. \eqref{eq weak form non hom} reaches the form
\begin{equation} \label{eq weak form}
- \langle
\dot{u} (\. ) , \dot{\f} (\. )
\rangle 
+ k 
\langle
u (\. ) , \f (\. )
\rangle =
\langle
f (\. ) , \f (\. )
\rangle
- \langle
\dot{v}_{0,\bar{T}}(\. )
, \dot{\f} (\. )
\rangle 
+ k 
\langle
v_{0,\bar{T}}(\. ) , \f (\. )
\rangle
\end{equation}

Finally, under the weak assumption that $e^{c\. }f(\. )$ is the weak derivative of a
function $F(\. ) \in L^{2}(I)$, assumption equivalent to $f \in H^{-1}$,
and by performing yet another integration by parts, one obtains
\begin{equation} \label{eq weak form estimates}
- \langle
\dot{u} (\. ) , \dot{\f} (\. )
\rangle _{c}
+ k 
\langle
u (\. ) , \f (\. )
\rangle _{c}=
- \langle
F (\. ) , \dot{\f}(\. )
\rangle _{0}
- \langle
\dot{v}_{0,\bar{T}}(\. )
, \dot{\f} (\. )
\rangle _{c}
+ k 
\langle
v_{0,\bar{T}}(\. ) , \f (\. )
\rangle _{c}
\end{equation}
The following fundamental lemma establishing the equivalence between $F(\. )$ and
$e^{c\.}f(\. )$ is relevant.
\begin{lemma} \label{lemequiv f and F}
Let $f \in H^{-1}$ and $F \in L^{2}$ be a week integral of $e^{c\.}f(\. )$. Then,
\begin{equation}
\|e^{c\.}f(\. ) \|_{H^{-1}} = \|F\|_{L^{2}} \leq K \|f \|_{H^{-1}}
\end{equation}
where $K>0$ is a constant depending only on $c$ and $\T$.
\end{lemma}
\begin{proof}
By definition and by integration by parts,
\begin{equation}
\begin{array}{r@{}l}
\|e^{c\.}f(\. ) \|_{H^{-1}}
&=
\sup \limits _{\substack{\f \in H^{1}_{0} \\ \|\f ' \|_{L^{2}} = 1}}
\int _{0}^{\T}e^{c\.}f(\. ) \f (\. )
\\
&=
\sup \limits _{\substack{\f \in H^{1}_{0} \\ \|\f ' \|_{L^{2}} = 1}}
\int _{0}^{\T}F(\. ) \f '(\. )
\\
&=\|F\|_{L^{2}}
\end{array}
\end{equation}
The second part follows easily since $\T < \infty$ and $e^{c\.}$ is a bounded smooth
function.
\end{proof}

The following lemma will be useful later on. Its proof is straightforward given the calculations
already carried out.
\begin{lemma} \label{lemcalc int force}
    Let $F = \int e^{c\.} \FF (x) $, i.e. $F$ be the primitive of the excitation function
    corresponding to the displacement function $x$ (multiplied by the corresponding exponential factor). Then,
    \begin{equation}
        F = F_{c,k}(x) = e^{c\.}\dx + k \int e^{c\.} x
    \end{equation}
    modulo a constant.
\end{lemma}

In practice, i.e. in the construction of numerical methods, the form of eq.
\eqref{eq weak form} is more useful, as it involves the excitation force function itself.
However, the form \eqref{eq weak form estimates} is more useful in the estimation of
error rates, since it places $f \in H^{-1}(I)$,\footnote{The space $H^{-1}(I)$ consists
of functions that are weak derivatives of functions in $L^{2}(I)$. The assumption that
$e^{c\. }f(\. )$ be in $H^{-1}$ is clearly equivalent to $f \in H^{-1}$, since
$e^{c\. }$ is bounded and smooth. The space $H^{-1}$ is the dual of $H^{1}$, i.e. the space
of all continuous linear functionals $H^{1} \ra \R$, see \cite{BrezisAnFonc}.} which is the
correct functional space for the excitation function, as it is the weakest regularity of an
excitation function for which the weak form and the Duhamel formula produce
a pointwise well-defined displacement function.

We can now state and prove the main theoretical result of the present paper, the theorem
establishing that eq. \eqref{eq weak form} is the weak formulation of the SDOF problem.

\begin{theorem} \label{thm weak form SDOF}
Eq. \eqref{eq weak form} is the \textit{weak form} of
eq. \eqref{eq SDOF}, with boundary conditions those of eq. \eqref{eq bound cond}. Equivalently
\begin{enumerate}
\item If $\bar{T}$ is fixed and does not belong to the
exceptional set $\{ \tau \in \R , \sin (\w _{d} \tau ) = 0 \}$, and if
$v_{0,\bar{T}} \in H^{1}$ satisfies the boundary conditions of eq. \eqref{eq bound cond}, then
for every $f \in H^{-1}$ the weak form admits a unique
solution $u \in H^{1}_{0}$
\item Moreover, if $f \in H^{s}(\bar{I})$, for some $s \geq -1$, then $u \in H^{s+2}$
\item the boundary value problems solved by the weak formulation, when we let the
function $v_{0,\bar{T}} \in C^{s+2}$ vary, are in a one-to-one correspondence
with the initial condition problems solved by the ODE of eq. \eqref{eq SDOF}.
\end{enumerate}
\end{theorem}

\begin{remark}
For the initial and final velocity to be well defined, the excitation function $f$ needs
only be $H^{-1/2 + \e }$ for any $\e >0$ locally around $0$ and $\T$, i.e. in some
intervals $[0,\d ]$ and $[\T -\d , \T]$ for some $\d >0$. This assumption allows for the
weak integral of $f(\.)$ to admit a trace operator, i.e. for its values on the boundary
of $I$ to be well defined.

This additional regularity of $f$ will be a standing assumption from now on. It is not
included in the notation, in order to keep it as light as possible. If $f$ is supposed
to be compactly supported, the assumption is automatically satisfied as soon as 
the interval $(0,\T)$ contains the support of $f$.
\end{remark}

The theorem states that eq. \eqref{eq weak form} is defined for less regular functions
$f$ (in our case for $f \in H^{-1}$) than the differential equation of eq. \eqref{eq SDOF}.
The two equations, eq. \eqref{eq weak form} and eq. \eqref{eq SDOF}, are equivalent
whenever they are both defined, in the sense that when $f$ is regular
enough ($f \in C^{0}$ will suffice), the solution of the weak form of the
equation is sufficiently regular in order to satisfy eq. \eqref{eq SDOF} and is
actually equal to its solution. Eq. \eqref{eq weak form} is indeed the weak form
of eq. \eqref{eq SDOF} since the two following facts hold true.
\begin{enumerate}
\item Firstly, the weak form is
equivalent to the classical one, under
the relevant regularity assumption on $f$.
\item Secondly only the first derivative
appears in the weak form, instead of the second derivative appearing in the ODE, since
one derivative has been transferred to the test function by integration by parts.
\end{enumerate}

The weak form is an equation whose unknown is the function $u \in H^{1}_{0}(I)$ and
the equation should be satisfied for every $\f \in C^{\infty}_{0}(I)$.
The function $v_{0,\bar{T}}$ is to be considered as a moving observer
whose position coincides with the moving body at times $0 $ and $\bar{T}$.
The unknown function $u$ is in fact the position of the body relative
to the moving observer.
The terms $\langle
\dot{v}_{0,\bar{T}}(\. )
, \dot{\f} (\. )
\rangle 
- k 
\langle
v_{0,\bar{T}}(\. ) , \f (\. )
\rangle$
on the right-hand-side of eq. \eqref{eq weak form} are to be considered as
the effects of fictitious forces in the (not necessarily inertial)
frame of reference attached to the moving observer.

The role of the functions $\f $ is that of test functions.
Evaluation of the integrals in eq. \eqref{eq weak form} can be
viewed as an act of measurement relative to the shape of the function
$\f$. On the right-hand-side, the correlation between the function $f$
and the test function is measured, by calculating
a virtual-work-like integral (it is actually the projection of the
function $f$ on the subspace generated by $\f$ in the sense
of the $L^{2}_{c}$ inner product). The weak form of the equation states
that this measurement should be equal to the one on the left-hand-side,
where the correlation between the test function and the position,
and the derivative of the test function and the velocity are measured.
If all such measurements agree, for all admissible test functions, then
the function $u$ solves the weak form of the equation.

It is also useful for later use to note that in the language of $L^{2}$ products,
$H^{-1}$ and $H^{1}$ duality, and of the operator $\FF$, the weak formulation amounts to
the equation
\[ \label{eqdef duaity}
\begin{array}{r@{}l}
_{H^{-1}}\langle e^{c\.}\FF (x)(\.),\f \rangle _{H^{1}} &=
\langle
f (\. ) - \FF (v_{0,\bar{T}})(\. ) , \f (\. )
\rangle _{c} 
\\
&=
\langle
e^{c\.}(f (\. ) - \FF (v_{0,\bar{T}})(\.)) , \f (\. )
\rangle _{0}
\\
&=
\langle
\tilde{F} , \dot{\f} (\. )
\rangle _{0} 
\end{array}
\]
where $\tilde{F}$ is an integral of $e^{c\.}(f (\. ) - \FF (v_{0,\bar{T}})(\.))$.

\begin{proof} [Proof of theorem \ref{thm weak form SDOF}] \renewcommand{\qedsymbol}{}
The existence of a solution follows from the Duhamel integral formula,
eq. \eqref{eq Duh}, which defines a (weak) solution $x \in H^{1}$ as long as
$f \in H^{-1}$. This is so because the Duhamel Integral representation, c.f. eq.
\eqref{eq Duh}, implies directly that
\begin{equation} \label{eqduh developed displ}
x (t)
=
\frac{e^{ - c t/2}}{{\omega _d }} \left(
\sin [\omega _d t]
\int_{0}^{t} e^{ c \tau} f (\tau )e^{ -c\tau /2} \cos [\omega _d \tau ]\dd\tau -
\cos [\omega _d t]
\int_{0}^{t} e^{ c \tau} f (\tau )e^{ -c\tau /2} \sin [\omega _d \tau ]\dd\tau
\right)
\end{equation}
which is in $L^{2}$ as soon as $f \in H^{-1}$ and integration is to be interpreted as
$_{H^{-1}}\langle \. , \. \rangle _{H^{1}}$ duality. The following lemma establishes this
fact. The calculations involved in its proof will be useful in \S \ref{secbounds theory}.

\begin{lemma} \label{lemDuhamel}
The solution is independent of the constant of integration when choosing the integral of
$e^{ c \. } f (\. )$ in the Duhamel integral formula.
\end{lemma}
\begin{proof}[Proof of the lemma]
Let $F(\. ) \in L^{2}_{0}$ be the weak integral of $e^{ c \. } f (\. ) \in H^{-1}$
satisfying $F(0) = 0$ and let $K \in \R$ be an arbitrary constant of integration. Call,
also,
\begin{equation}
\begin{array}{r@{}l}
d(t) &= e^{ -ct /2} \cos [\omega _d t ], d(0) = 1, \dot{d}(0 ) = -c/2 \\
s(t) &= \frac{e^{ -ct /2}}{\w _{d}} \sin [\omega _d t ] \phantom{]}
, s(0) = 0, \dot{s}(0 ) = 1
\end{array}
\end{equation}
the two fundamental solutions of the SDOF problem. Direct calculation
yields
\begin{equation} \label{eqder fund sol}
\begin{array}{r@{}l}
\dot{d}(t) &= -\frac{c}{2}d(t) - \w _{d}^{2} s(t) \\
\dot{s}(t) &= -\frac{c}{2}s(t) + d(t) 
\end{array}
\end{equation}
With this notation, eq. \eqref{eqduh developed displ} gives
\begin{equation} \label{eqsol Duh F}
\begin{array}{r@{}l}
x(t) &=
s(t)\int_{0}^{t} e^{ c \tau  } f (\tau  ) d(\tau ) \dd\tau
-
d(t)\int_{0}^{t} e^{ c \tau  } f (\tau  ) s(\tau ) \dd\tau
\\
&=
s(t)\int_{0}^{t} \frac{\mathrm{d}}{\mathrm{d}t}(F(\tau  )+ K) d(\tau ) \dd\tau
-
d(t)\int_{0}^{t} \frac{\mathrm{d}}{\mathrm{d}t}(F(\tau  )+ K) s(\tau ) \dd\tau
\\
&=
s(t)
\left(
(F(t)+K)d(t)- K - \int_{0}^{t} (F(\tau  )+ K) \dot{d}(\tau ) \dd\tau
\right) - \\
&\phantom{s(t)\int_{0}^{t}s(t)\int_{0}^{t} } -
d(t)
\left(
(F(t)+K)s(t) - \int_{0}^{t} (F(\tau  )+ K) \dot{s}(\tau ) \dd\tau
\right)
\\
&=
-s(t)\int_{0}^{t} F(\tau  ) \dot{d}(\tau ) \dd\tau
+ 
d(t)
 \int_{0}^{t} F(\tau  ) \dot{s}(\tau ) \dd\tau
\\
&=
s(t)
\int_{0}^{t} F(\tau  ) (\frac{c}{2}d(\tau) + \w _{d}^{2} s(\tau) ) \dd\tau
+
d(t)
 \int_{0}^{t} F(\tau  ) (d(\tau) -\frac{c}{2}s(\tau) ) \dd\tau
\end{array}
\end{equation}

This ends the proof of the lemma, since the last three expressions are independent
of $K$. The formula in the last expression will be useful in \S \ref{secbounds theory}.
\end{proof}

Back to the proof of the theorem, one can obtain the velocity $\dx$ by
taking one derivative of the Duhamel integral. The same type of calculation gives
\begin{equation} \label{eqduh developed speed}
\begin{array}{r@{}l}
\dx(t) &= \frac{\dd}{\dd t} \left(
s(t)\int_{0}^{t} e^{ c \tau  } f (\tau  ) d(\tau ) \dd\tau
-
d(t)\int_{0}^{t} e^{ c \tau  } f (\tau  ) s(\tau ) \dd\tau
\right)
\\
&=
e^{-ct}F(t) + (d(t) -\frac{c}{2}s(t) )
\int_{0}^{t} F(\tau  ) (\frac{c}{2}d(\tau) + \w _{d}^{2} s(\tau) ) \dd\tau
\\
&\phantom{(d(t) -\frac{c}{2}s(t) )
\int_{0}^{t} F(\tau  ) (\frac{c}{2}d(\tau)}
-(\frac{c}{2}d(t) + \w _{d}^{2} s(t)
 \int_{0}^{t} F(\tau  ) (d(\tau) -\frac{c}{2}s(\tau) ) \dd\tau 
\end{array}
\end{equation}
which again is in $L^{2}$ as soon as $f \in H^{-1}$, which directly yields that $x \in \HH$.

Calculation of the acceleration in the same fashion cannot result in anything better
than
\begin{equation} \label{eqduh developed acc}
\ddot{x} (t) = 
 f (t ) - x(t) - c\dx (t)
\end{equation}
which has the same regularity as $f$, i.e. is in $H^{-1}$.

\bigskip

The equivalence of the weak and the strong forms of the SDOF problem follows from
eq. \eqref{eq weak form estimates} which implies that
\begin{equation}
\langle
\dot{u} (\. ) , \dot{\f} (\. )
\rangle _{c}=
 k 
\langle
u (\. ) , \f (\. )
\rangle _{c} +
\langle
F (\. ) , \dot{\f} (\. )
\rangle _{0}
- \langle
\dot{v}_{0,\bar{T}}(\. )
, \dot{\f} (\. )
\rangle _{c}
+ k 
\langle
v_{0,\bar{T}}(\. ) , \f (\. )
\rangle _{c}
\end{equation}
where $u \in H^{1}_{0}$, and $v_{0,\bar{T}}(\. ) $ is smooth by construction.
If $f \in H^{-1}$ or equivalently $F \in L^{2}$, then the above equation implies that
$\dot{u} \in L^{2}$, and therefore $u \in H^{1}_{0}$. By
induction, if $f \in H^{s}$, $s \geq -1$, then $u \in H^{s+2} \cap H^{1}_{0}$ and the
transformation leading from eq. \eqref{eq SDOF} to the eq. \eqref{eq weak form}
can be reversed, showing that $x(\. ) = u(\. ) + v_{0,\bar{T}}(\. )$
is a solution to the initial value problem.

Finally, uniqueness for the weak form reduces to uniqueness for the classical problem, thanks to
the classical argument for uniqueness of solution of linear systems. Let
$u_{i}$, $i=1,2$ be a solution to eq. \eqref{eq weak form} for the same excitation
function $f$, initial conditions and final time $\bar{T}$, and therefore also for the same
boundary conditions. Calling $u = u_{1} - u_{2} \in H^{1}_{0}$, one obtains directly
\begin{equation}
- \langle
\dot{u} (\. ) , \dot{\f} (\. )
\rangle 
+ k 
\langle
u (\. ) , \f (\. )
\rangle =
0
\end{equation}
for all test functions $\f$. Consequently,
\begin{equation} \label{eq. b}
- \langle
\dot{u} (\. ) , \dot{\f} (\. )
\rangle 
=- k 
\langle
u (\. ) , \f (\. )
\rangle
\end{equation}
which implies that $\dot{u}$ has a weak derivative which can be calculated by integration
by parts, yielding
\begin{equation}
\ddot{u} (\. )= -ku (\. ) - c\dot{u}(\. )
\end{equation}
coupled with $0 $ boundary conditions: $u(0) = u(\bar{T}) = 0$. Since the right-hand side of eq.
\eqref{eq. b} is continuous, the weak derivative $\ddot{u}$ is in fact a strong derivative
and by induction $u \in C^{\infty}$. The choice of $\bar{T }$ now implies that $u 
\equiv 0$, which is equivalent to the uniqueness of the weak solution.
\end{proof}

\begin{remark}
The test functions need only be of regularity $H^{s+2}_{0}$ when
$f \in H^{s}$.
\end{remark}
This remark is important for numerical applications, since no
additional smoothness condition needs to be imposed, other than, say,
continuity of the first derivative.

\section{Covariance of the Weak form with respect to time reparametrization} \label{sec scaling}

If one considers the SDOF equation in its classical form with unit mass and restricted
on a finite interval, cf. eq. \eqref{eq SDOF}, six defining parameters appear in the
problem: $c,k,f,x_{0},\dx_{0}$, and $\T$, the mass having been eliminated via
normalization.

Upon introduction of time-dilation $t \mapsto \l ^{-1} t$, one more parameter can be
eliminated, and the remaining parameters are transformed as follows. Introducing the
local notation $\chi (t) = x(\l ^{-1}t)$ and $\f (t) = f(\l ^{-1}t)$, one gets that
\begin{equation}
\l^{2} \ddot{\chi}(\l t)  + \l c \dot{\chi}(\l t) + k \chi(\l t)  = f(t) \text{ in } [0,\bar{T}]
\end{equation}
or, equivalently,
\begin{equation}
\ddot{\chi}  + \l^{-1} c \dot{\chi} + k \l^{-2} \chi  = \l^{-2} \f \text{ in } [0,\l ^{-1}\bar{T}]
\end{equation}
Under time-dilation, therefore, the parameters defining the SDOF problem are transformed
by
\begin{equation} \label{eqtime reparametrization}
(c,k,f,x_{0},\dx_{0},\T)
\mapsto
(\frac{c}{\l}, \frac{k}{\l^{2}},
\frac{f(\l^{-1}\.)}{\l^{2}},x_{0},
\frac{\dx_{0}}{\l}, \frac{\T }{\l} )
\end{equation}
In particular, for $\l = \w_{n}= \sqrt{k}$, one gets that $k$ is normalized to $1$ and
$c$ to $c/\sqrt{k} = 2\xi$, cf. eqs. \eqref{eq damp} and \eqref{eq undamp eigenfr}.

It can be easily verified that the weak formulation transforms in the same way 
under time-dilations. Therefore, the study of the weak formulation can be restricted
to $k=1$, provided that all other parameters are not constrained. This will be helpful
in obtaining error estimates for methods constructed upon the weak formulation.

\section{Relation with Lagrangian and Hamiltonian mechanics}

For this section see e.g. \cite{AbrahamMechanics}. See, also,
\cite{Lager88}, section 2.2, for a similar manipulation.

The Lagrangian for the undamped linear oscillator with $0$ excitation function is
defined as
\begin{equation} 
\begin{array}{r@{}l}
\LL (x,\dot{x}) = \LL_{0} (x,\dot{x}) &= \frac{1}{2} \dot{x}^{2}  -
\frac{1}{2} kx ^{2}\\
&= T(\dot{x}) - V(x)\\
&= T_{0}(\dot{x}) - V_{0}(x)
\end{array}
\end{equation}
where $T $ is the kinetic and $V$ is the potential energy.
The equation of motion can be obtained by the formula
\begin{equation} \label{eq mot lag}
\frac{\dd}{\dd t}\frac{\partial}{\partial \dot{x}}\LL (x,\dot{x}) =
\frac{\partial}{\partial x}\LL (x,\dot{x})
\end{equation}
and is coupled with boundary conditions as in eq. \eqref{eq bound cond}. Solutions to
the weak
problem for the homogeneous system are obtained as the stationary points of the averaged
action functional
\begin{equation}
\int _{0}^{\T }\LL (x,\dot{x})
\end{equation}
The differential of this functional can be easily calculated as
\begin{equation}
\mathrm{D}\left(
\int _{0}^{\T }\LL (x,\dot{x}) 
\right)
\. \f = \langle \dx , \dot{\f} \rangle_{0}
-
\langle x , \f \rangle_{0} , \forall \f \in \HH _{0}
\end{equation}
Stationary points are therefore given by the equation
\begin{equation}
\mathrm{D}\left(
\int _{0}^{\T }\LL (x,\dot{x}) 
\right)
\. \f = 0 , \forall \f \in \HH _{0}
\end{equation}
and an integration by parts with respect to $\dx$ gives the equation of motion, i.e.
the strong SDOF problem.

The Hamiltonian function for the same system is
\begin{equation} \label{ptwise hamiltonian unrect}
\begin{array}{r@{}l}
\Ham (x,\dot{x}) = \Ham _{0}(x,\dot{x}) &= \frac{1}{2} \dot{x}^{2} +
\frac{1}{2} kx ^{2}\\
&= T(\dot{x}) + V(x)
\end{array}
\end{equation}
The equation of motion can be obtained by the formula
\begin{equation}
\begin{array}{r@{}l}
\frac{dx}{dt}  &= \phantom{-} \frac{\partial \mathcal{H} (x,\dx)}{\partial \dx}\\
\frac{d \dx}{dt}  &= -\frac{\partial \mathcal{H} (x,\dot{x})}{\partial x}
\end{array}
\end{equation}
and is coupled with initial conditions as in eq. \eqref{eq in cond}.

Both formulations work for the undamped system where total energy is conserved, which in
this notation is equivalent to saying that
\begin{equation}
\frac{\dd}{\dd t}\mathcal{H} (x,\dx) = 0
\end{equation}
provided that $(x,\dx)$ is a solution of the equation of motion, eq. \eqref{eq SDOF}.
Direct calculation can show the more general fact that solutions of the Lagrangian
problem have constant energy.

In the context of a damped system, where energy is dissipated due to the presence
of the term $c \dx ^{2}$ in the relevant calculations (see \S \ref{seccdis law}),
and in the presence of an excitation function, where energy input takes place, the
Lagrangian and Hamiltonian formulations can be modified in the following way.

The time-dependent Lagrangian for the free SDOF system is
defined as
\begin{equation} \label{ptwise lagrangian}
\begin{array}{r@{}l}
\LL _{c} (x,\dot{x},t ) &= \frac{1}{2} e^{ct} \dot{x}^{2} -
\frac{1}{2} ke^{ct}x ^{2}\\
&= T_{c}(\dx ,t) - V_{c}(x ,t)
\end{array}
\end{equation}
The equation of motion for the free damped system can be obtained by the formula of eq.
\eqref{eq mot lag}, coupled with boundary conditions.

The weak formulation for a system with excitation function $f$ can be obtained in
otherwise the same way by solving the equation
\begin{equation}
\mathrm{D}\left(
\int _{0}^{\T }\LL _{c} (x,\dot{x}) 
\right)
\. \f = \langle f , \f \rangle _{c}
= \langle F , \dot{\f} \rangle _{0} , \forall \f \in \HH _{0}
\end{equation}
Instead of looking for stationary points, one asks, thus, that the differential be
equal to a certain functional in $L^{2}_{c}$. The equation of motion is once again
obtained by integration by parts.

In other words, the function $x$ is a solution to
the weak problem if, and only if, for fixed boundary conditions, the
variation of the functional $\LL $ along any virtual deformation
$\e \f $ of the orbit $x$ is equal to the virtual work
$\e \langle F  , \dot{\f}  \rangle _{0}$, up to terms
of the order of $\e ^{2}$. More precisely, the condition is that,
for every $\f \in H^{1}_{0}$,
\begin{equation}
\LL _{c} (x+ \e \f ) =\LL _{c} (x ) + \e \langle F  , \dot{\f}  \rangle _{0}
+ O(\e ^{2})
\end{equation}

The time-dependent Hamiltonian function for the same system is
\begin{equation} \label{ptwise hamiltonian}
\begin{array}{r@{}l}
\mathcal{H}_{c} (x,\dot{x}) &= \frac{1}{2} e^{ct} \dot{x}^{2} (t) +
\frac{1}{2} e^{ct} kx ^{2}(t)\\
&= T_{c}(\dot{x}) + V_{c}(x)
\end{array}
\end{equation}
The equation of motion for the free system can be obtained by the formula of eq.
\eqref{ptwise hamiltonian unrect}, as can show a simple calculation.

In the presence of dissipation, the time-dependent Hamiltonian function of the
homogeneous problem is not preserved, but satisfies the equation
\begin{equation}
\frac{\dd}{\dd t}\Ham _{c} (x,\dx,t) = - c \LL _{c} (x,\dx,t) %
\end{equation}
as can show a direct calculation. Therefore, the role of the exponential factor
$e^{ct}$, introduced in eq. \eqref{ptwise lagrangian}, is to balance the dissipation of
energy due to damping, since the total energy $\Ham _{c}$ is expected to oscillate
around its initial value, while $\Ham _{0} \ra 0$ as $t \ra \infty$ for the same system,
since $\Ham _{0}$ represents the standard mechanical energy and the system dissipates
energy.

\subsection{A conservation law} \label{seccons law}

In the case where $u \in H^{2} \cap H^{1}_{0}$, the following proposition can be proved.
It is generalization of the conservation of energy in the dissipative setting, since in the
case where $c = 0$, the energy conservation law is obtained.
Here, no hypothesis on $\T $ is necessary.

\begin{prop}
Let $f \in L^{2}$and $u  \in H^{2} \cap H^{1}_{0}$ be the solution of the weak SDOF problem
with $0$ boundary conditions. Then,
\begin{equation}
\frac{1}{2}e^{c\T }\dot{u}(\T )^{2} - \frac{1}{2}\dot{u}(0 )^{2}
=
\langle f , \dot{ u} \rangle _{c} - 
c \int _{0}^{\T } \LL _{c} (u,\dot{u})
= \langle \dot{ f} , u \rangle _{c}
\end{equation}
\end{prop}
In the case where $c=0$, the conclusion of the lemma reads simply
\begin{equation}
\begin{array}{r@{}l}
\frac{1}{2}\dot{u}(\T )^{2} - \frac{1}{2}\dot{u}(0 )^{2}
&=
\int _{0}^{\T } f  \dot{ u}\\
&= \int _{u(\. )} f  du
\end{array}
\end{equation}
which is precisely energy conservation in view of the hypothesis $u(0) = u(\T ) = 0$.

In the dissipative case, $c>0$, the conclusion of the proposition can be obtained
via the Lagrangian and Hamiltonian formulations.

Since one can transform non-homogeneous boundary conditions to homogeneous one
by modifying the excitation function $f$, the corresponding hypothesis poses no real
restriction to applications.

\begin{proof}
Multiply the equation
\begin{equation}
d(e^{ct }\dot{u}(t)) + e^{ct } u(t) = e^{ct }f(t)
\end{equation}
by $\dot{u}(t) \in H^{1}$ and integrate by parts. This gives
\begin{equation}
\begin{array}{r@{}l}
e^{ct }\dot{u}(t)^{2} \rvert_{0}^{\T } 
-
\langle \dot{ u} , \ddot{ u} \rangle _{c}
-
c k \| u\|_{c}^{2}
-
k \langle u , \dot{ u} \rangle _{c}
&=
\langle f , \dot{ u} \rangle _{c} \iff \\
e^{ct }\dot{u}(t)^{2} \rvert_{0}^{\T } 
-
c k \| u\|_{c}^{2}
-
\langle \dot{ u} , \ddot{ u} + ku \rangle _{c}
&=
\langle f , \dot{ u} \rangle _{c}\iff \\
e^{ct }\dot{u}(t)^{2} \rvert_{0}^{\T } 
-
c k \| u\|_{c}^{2}
-
\langle \dot{ u} , f - c\dot{ u} \rangle _{c}
&=
\langle f , \dot{ u} \rangle _{c}\iff \\
e^{ct }\dot{u}(t)^{2} \rvert_{0}^{\T } 
-
2c \int _{0}^{\T }\LL _{c} (u,\dot{u})
&=
2\langle f , \dot{ u} \rangle _{c}
\end{array}
\end{equation}
Where use has been made of the fact that eq. \eqref{eq SDOF} is satisfied in the weak
sense. This gives directly the first part of the equality to be proved.

The second part follows from another integration by parts applied on the term
$\langle f , \dot{ u} \rangle _{c}$, and from the fact that
$\LL  _{c}(u,\dot{u}) = \langle f , u \rangle _{c}$.
The details are left to the reader.
\end{proof}


Dropping the homogeneous boundary condition hypothesis, one can obtain
the following corollary, practically a restatement of the one here
above.
\begin{corollary}
Let $x \in H^{2}$ and $f \in \HH$ be the solution and the excitation function of the
weak formulation of the SDOF problem. Then,
\begin{equation}
\begin{array}{r@{}l}
| \Ham _{c} (\T ) - \Ham _{c} (0)| &\leq \langle f , \dot{ x} \rangle _{c}
-
c\langle f , x \rangle _{c}
\\
&\leq
e^{c\T }f(\T )x(\T ) - f(0)x(0) -  \langle \dot{ f} , x \rangle _{c}
\end{array}
\end{equation}
\end{corollary}

In the particular case of $0$ initial conditions, the above corollary yields the following
inequality:
\begin{equation}
| \Ham _{c} (\T )| \leq
e^{c\T }f(\T )x(\T ) -  \langle \dot{ f} , x \rangle _{c}
\end{equation}

\subsection{A dissipation law} \label{seccdis law}

The strong form of the SDOF problem can be rewritten as
\begin{equation}
\ddot{x} + k x = f - c \dx
\end{equation}
Multiplication by $\dx$ and integration by parts yields the total loss of
energy in time $t \in I$ (cf. eq. \eqref{ptwise hamiltonian} for the definition of
$\Ham = \Ham _{0} $):
\begin{equation}
\begin{array}{r@{}l}
\frac{1}{2}\dx ^{2} \rvert _{0}^{t} + \frac{1}{2} k x ^{2} \rvert _{0}^{t}
&=
\int _{0}^{t } f \dx - c \int _{0}^{t } \dx ^{2} \iff \\
\Ham (t) - \Ham (0)
&=
f(t)x(t) - f(0)x(0) - \int _{0}^{t } \dot{f} x - c \int _{0}^{t } \dx ^{2}
\end{array}
\end{equation}

The following proposition follows immediately from this dissipation law.
\begin{prop}
If $x \in H^{2}$ and $f \in L^{2}$ are respectively the solution and the excitation function
of a weak SDOF problem. Then
\begin{equation}
\frac{1}{2}\| \dx \|_{c}^{2} + \frac{1}{2} k \|x\|_{c}^{2} \leq
\left(
\|f\|_{c} + \frac{1}{c}(e^{c\T }-1)\|\dot{f}\|_{c}
\right)
\|x\|_{c} +
\frac{1}{c}(e^{c\T }-1)(\Ham (0) -f(0)x(0))
\end{equation}
\end{prop}
\begin{proof}
Suppose that $c > 0$. The dissipation law gives directly
\begin{equation}
\begin{array}{r@{}l}
\int _{0}^{\T } e^{ct} \Ham  (t)
&\leq
\int _{0}^{\T }e^{ct} f(t)x(t) -
\int _{0}^{\T }e^{ct} \left(
\int _{0}^{t } \dot{f}(t)x(t)
\right)
+\frac{1}{c}(e^{c\T }-1)(\Ham (0) -f(0)x(0))
\\
&\leq
\left(
\|f\|_{c} + \frac{1}{c}(e^{c\T }-1)\|e^{-c\. }\dot{f}(\. )\|_{c}
\right)
\|x\|_{c} +
\frac{1}{c}(e^{c\T }-1)(\Ham (0) -f(0)x(0))
\end{array}
\end{equation}
where $\frac{1}{c}(e^{c\T }-1)$ is of course the measure of the interval $I$ with respect
to the measure $e^{ct}dt$ when $c>0$ and the H\"{o}lder inequality has been applied
twice. One can formally pass to the limit as $c \ra 0$ and obtain the factor
$\T$ in the case $c=0$.
\end{proof}

%
%

\section{Some remarks on the weak formulation}

The weak formulation has an inherent non-local character.
The property that characterizes the solution is a non-local one,
obtained by integration of the solution against test functions, and
therefore the stepwise character of the traditional methods is lost.
The loss of the local character is underlined by the fact that the
weak formulation works only with a boundary value problem, so that
the final displacement of the moving body is considered as a known
quantity. A function in $H^{1}$ does not have a well-defined pointwise
derivative, so the initial condition $\dot{x}_{0}$ makes
no sense in this context. Finally, no discretization in the time
domain is necessary for obtaining the weak formulation, but it can
be eventually introduced a posteriori for applications.

What is traditionally called weak formulation (see, e.g., \cite{Hugh87})
is in fact the projection of the (full) weak formulation, as it is
obtained in the present work, on a finite
dimensional space via the following procedure. A finite dimensional
subspace of $\HH $ is constructed by introducing a discretization of
the time domain, and by choosing the order of polynomial approximation,
typically equal to $2$. Then, the solution is \textit{a priori}
projected on this finite dimensional space of functions that are
$C^{2}$ and piece-wise polynomial of the given order, with
breakpoints at the points of the time-discretization mesh. The
dimension of the function space thus constructed can be computed, but
it is irrelevant, the important point being that it is of positive
(in fact infinite)
codimension\footnote{The codimension of a subspace of a
vector space is the dimension of its complementary space. In this
context, the point of the argument is that, since the codimension of
the space on which the solution is projected is infinite, the
information lost a priori in the construction of the weak solution
lives in an infinite dimensional space.}. The construction of the
projection, which to the best knowledge of the authors seems to be
arbitrary, consists essentially in choosing and fixing a weight
function (typically denoted by $w$). This weight function replaces the
test function $\f $ of this work and it is \textit{fixed} for each
numerical method, as its choice actually defines each particular
method. However, in
this work the test function cannot be fixed but has to
run through the entire function space of admissible weights.
The choice of the weight function determines the scope and the
efficiency of the method, under the limitations related to the
infinite codimensionality of the space in which the solution is
sought and the limitation in the choice of the (unique) test function.
These limitations are exactly those described by the Dahlquist Barrier
Theorem, (see \cite{Dahl56} and \cite{Dahl63}).

In this work, the actual weak formulation is obtained, since the
eq. \eqref{eq weak form} is a representation of the solution in the
full space $\HH $, and the construction of a projection onto a finite
dimensional subspace is left to be introduced a posteriori for the
eventual numerical implementations. Here, the projection cannot be
chosen arbitrarily, since it is dictated by the $L^{2}$ scalar product
of eq. \eqref{eq def in prod c}.
The arbitrary choice defining the particular numerical method would
consist in the choice of the space onto which the solution is
projected. A reasonable choice seems to be the space of $C^{1}$,
piece-wise polynomial functions of constant but arbitrary order.
Such a method is constructed in part two of the paper.

Finally, the weak formulation has the advantage that only one derivative
appears, and the second derivative is replaced by integration of the
unknown function against test functions. Therefore, the price to pay
for dropping the degree of the equation by one is the infinite
dimensionality of the space in which the problem is posed, the space
$\HH$. On the
other hand, this very fact allows higher-degree polynomial
approximation, and derivation is replaced by integration, which is
numerically more stable.

\section{Weak formulation: The MDOF case}
The exact same reasoning applies to the Multi-Degree-Of-Freedom (MDOF)
case, where the governing equation is
\begin{equation} \label{eq MDOF}
[M] \{ \ddot{x} \} + [C] \{ \dot{x} \} + K \{ x \}  = \{ f \} 
\end{equation}
where $[M],[C],[K] \in M_{n}(\R )$ are symmetric real $n \times n$ matrices.
The mass matrix $[M] $ is symmetric and positive definite, as is the
stiffness matrix $[K]$. The unknown is now a vector-valued function
$\{ x \}(\. ) \colon \R \ra \R ^{n}$, where $n \in \N ^{*}$ is the number
of degrees of freedom of the system.

Normalization with respect to the mass in this context amounts to left
multiplication by the matrix $[M]^{-1}$ which gives
\begin{equation} \label{eq MDOF norm}
\{ \ddot{x} \} + [\bar{C} ] \{ \dot{x} \}  + [\bar{K} ] \{ x \}  = \{ \bar{f} \} 
\end{equation}
where
\begin{align}
[\bar{C} ] = [M]^{-1}[C] \\
[ \bar{K} ] = [M]^{-1}[K] \\
\{ \bar{f} \}= [M]^{-1}\{ f \}
\end{align}
Multiplication by $e^{[\bar{C} ]t}$ and use of the Leibniz rule gives
\begin{equation}
\frac{\dd}{\dd t}\left[e^{[\bar{C} ]t}\{ \xx \} (t) \right] +
 e^{[\bar{C} ]t} [\bar{K} ] \{ x \}(t) = e^{[\bar{C} ]t}\{ f \}(t)
\end{equation}
As for the SDOF case, left multiplication by $\{ \f \} ^{\mathrm{T}} $, where
$\{ \f \}  \in C^{\infty}_{0}(I,\R ^{n})$, and integration by parts yields
\begin{multline} \label{weak form mdof in cond}
-\int _{0}^{\bar{T}} \{ \dot{\f} \}^{\mathrm{T}} (t) e^{[\bar{C} ]t} \{ \xx \} (t)  \dd t +
 \int _{0}^{\bar{T}} \{ \f \} ^{\mathrm{T}} (t) e^{[\bar{C} ]t}[\bar{K} ]\{ x \}(t)  \dd t = \\
\int _{0}^{\bar{T}} \{ \f \} ^{\mathrm{T}} (t) e^{[\bar{C} ]t}\{ f \}(t)  \dd t
\end{multline}
Eq. \eqref{weak form mdof in cond} should be satisfied for every
$\{ \f \} \in C^{\infty}_{0}(I,\R ^{n})$. The compact form of eq. \eqref{weak form mdof in cond} is obtained by introducing the space
$L^{2}_{c}(I ,\R ^{n})$, consisting of the square-integrable
vector-valued functions on $I$ and its inner product
\begin{equation} \label{inner prof mdof}
\langle
\{ \f \} (\. ) , \{ \psi \} (\. )
\rangle =
\langle
\{ \f \} (\. ) , \{ \psi \} (\. )
\rangle _{[\bar{C} ]} =
\int _{0}^{\bar{T}}  \{ \f \} ^{\mathrm{T}} (t) e^{[\bar{C} ] t} \{ \psi \}  (t) \dd t 
\end{equation}
In view of eq. \eqref{inner prof mdof}, eq.
\eqref{weak form mdof in cond} can be rewritten as
\begin{multline} \label{eq weak form mdof non hom}
- \langle
\{ \xx \} (\. ) , \{ \dot{\f} \} (\. )
\rangle 
+ 
\langle
[\bar{K} ] \{ x \} (\. ) , \{ \f \} (\. )
\rangle =
\langle
\{ f \} (\. ) , \{ \f \} (\. )
\rangle ,  \\
\forall \{ \f \} \in C^{\infty}_{0}(I ,\R ^{n})
\end{multline}

Working in eigenmode coordinates $\phi _{i}$, with $i =
1, \cdots , n$, where the equations assume the decoupled
form
\begin{equation}
\mu _{i}\ddot{\phi}_{i}  + \xi _{i} \dot{\phi}_{i}  + \kappa _{i} \phi _{i}  = \bar{f}_{i} \text{ for } i = 1 , \cdots , n
\end{equation}
with initial conditions $\phi _{i} (0) = \phi _{i,0}$ and
$\dot{\phi}_{i}(0)= \dot{\phi }_{i,0}$, one can directly verify the
following facts. To each eigenmode $\phi _{i}$ corresponds, as in
eq. \eqref{eq bound prob}, an affine linear mapping
\begin{equation} \label{eq map bound val MDOF}
\BB _{\phi _{i,0},\bar{f} _{i},\bar{T}} \colon \dot{\phi}_{i,0} \mapsto x(\bar{T})
\end{equation}
Therefore, the initial conditions in physical coordinates,
$\{ x _{i,0} \}$ and $\{ \dot{x} _{i,0} \}$ can be transformed into
initial conditions in eigenmode coordinates
$\{ \phi _{i,0} \}$ and $\{ \dot{\phi } _{i,0} \}$, then to
boundary conditions $\{ \phi _{i,0} \}$ and $\{ \phi _{i,\bar{T}} \}$
and finally to boundary conditions in physical coordinates
$\{ x _{i,0} \}$ and $\{ x _{i,\bar{T}} \}$.
As in the SDOF case, each operator $\BB _{\phi _{i,0},\bar{f} _{i},\bar{T}}$
is invertible iff $\bar{T}$ is not a half-integer multiple of the
eigenperiod $T_{i}$ of the corresponding SDOF system. Consequently,
the correspondence between initial conditions and boundary value
for the MDOF problem (both in physical and eigenmode coordinates)
is bijective (i.e. one-to-one and onto) provided that
\begin{equation}
\bar{T} \notin \bigcup _{i=1}^{n} \frac{1}{2} \N ^{*} T_{i}
\end{equation}
This condition is again typical in any relevant sense, just as
for the SDOF problem.

Assuming given boundary conditions, one can fix
two $C^{\infty}$ matrix valued mappings $[ v _{0} ], [ v _{\bar{T}} ]\colon [0,\bar{T}] \ra M_{n}(\R) $ such that
\begin{equation}
\begin{cases}
v_{ii,0}(0) = v_{ii,\bar{T}}(\bar{T}) = 1, 1 \leq i \leq n \\
v_{ii,0}(\bar{T}) = v_{ii,\bar{T}}(0) = 0, 1 \leq i \leq n \\
v_{ij} (\. ) \equiv 0, i \neq j
\end{cases}
\end{equation}
Then, if $\{ x \}(\. ) $ solves eq. \eqref{eq MDOF} with boundary conditions
\eqref{eq bound cond}, then the function
\begin{equation}
 \{ u \}(\. ) = \{ x \}(\. ) -  [v_{0}  ](\. )\{ x \}_{0} -  [v_{\bar{T}} ](\. )\{ x \}_{\bar{T}}
\end{equation}
satisfies homogeneous boundary conditions $u(0)= u(\bar{T}) =0$.

Introduction of the mapping
\begin{equation}
\{ u \}(\. ) = \{ x \}(\. ) -  [v_{0}  ]\{ x \}_{0}(\. ) - [v_{\bar{T}} ] \{ x \}_{\bar{T}}(\. )
\end{equation}
yields the decomposition
\begin{equation}
 \{ x \}(\. ) = \{ u \} (\. )  + \{v_{0, \bar{T}} \}(\. )
\end{equation}
where, as for the SDOF case,
\begin{equation}
\{v_{0, \bar{T}} \}(\. ) =  [v_{0}  ](\. ) \{ x \}_{0} +  [v_{\bar{T}} ](\. ) \{ x \}_{\bar{T}}
\end{equation}
Then, eq. \eqref{eq weak form non hom} reaches the form
\begin{multline} \label{eq weak form mdof}
- \langle
\{ \dot{u} \} (\. ) , \{ \dot{\f} \} (\. )
\rangle 
+ 
\langle
[\bar{K} ] \{ u \} (\. ) , \{ \f \} (\. )
\rangle = \\
\langle
\{ f \} (\. ) , \{ \f \} (\. )
\rangle
+ \langle
\{ \dot{v}_{0,\bar{T}} \}(\. )
, \{ \dot{\f} \} (\. )
\rangle 
- 
\langle
[\bar{K} ] \{v_{0, \bar{T}} \}(\. ) , \{ \f \} (\. )
\rangle
\end{multline}
which should hold for all $\f \in C^{\infty}_{0}$.
\begin{theorem} \label{thm weak form MDOF}
Eq. \eqref{eq weak form mdof} is the \textit{weak form} of
eq. \eqref{eq MDOF}, with boundary conditions
$\{ x \} (0) =\{ x_{i,0} \}$ and
$\{ x \}(\bar{T}) =\{ x_{i,\bar{T}} \}$.
\end{theorem}

This is the weak form of the MDOF problem, and as in the SDOF case,
it is an
equation with unknown $\{ u \} (\. ) \in H^{1}_{0}(I , \R ^{n})$, and
should be satisfied for all $\{ \f \} (\. ) \in C^{\infty}_{0}
(I , \R ^{n})$.

The weak form of the MDOF problem has the obvious disadvantage that
it requires the calculation of the exponential matrix
\begin{equation}
e ^{[\bar{C } ]t }
\end{equation}
which is computationally efficient only if the matrix $[\bar{C} ] $ is
diagonal. Therefore, for numerical applications the reader should bear
in mind the implicit condition that the matrix $[\bar{C} ] = [M]^{-1}[C]$
is considered to be diagonal.

The authors plan to examine the case where $[\bar{C} ] $ is diagonal dominant
in a future article.


%
\section{From the weak formulation to numerical implementations: the SDOF case} \label{sec theory to alg}

An algorithm based on the weak formulation can be constructed by restricting the solution
to $\EH$, a finite dimensional subspace of $H_{0} ^{1} $. The restriction is obtained
via a projection of the force with respect to the $H^{-1}$ duality, and then obtaining
the solution associated to this problem.
The description that follows applies directly to the SDOF problem, and with some obvious
modifications to the MDOF one.

\subsection{Obtaining a numerical method}\label{seccons app sol}

Suppose that $\EH$ is generated by the $n $ functions
$\{ b_{i} \}_{i=1}^{n} \subset \scIh $, i.e.
\begin{equation}
\EH = \{ \f \in \HH _{0} , \exists  \b _{i} \in \R , \f (\. ) = \sum _{i=1}^{n}
\b _{i} b _{i} (\. ) \}
\end{equation}
For obvious reasons, it is assumed that there exists a function $b$ in $\EH$ such that
$\dot{b}(0) \neq 0$.

The unknown function $u $ is restricted to $\EH$, and given the choice
of a basis of $\EH$, namely $\{ b_{i} \}_{i=1}^{n}$, the approximate
solution $u_{ap}$ has an expression $u_{ap}(\. ) = \sum _{1}^{n} u_{i}b_{i} (\. )$.
The basis should be completed by two functions $b_{0}$ and $b_{n+1}$ such that
\begin{equation}
b_{0} ( 0) = b_{1} ( 1) = 1 \text{ and } b_{0} ( 1) = b_{1} ( 0) = 0
\end{equation}
(notice the abuse of notation since $b_{0}$ and $b_{1}$ are elements of $H^{1} \setminus H_{0}^{1}$).
These functions will be used to form $v_{0,\bar{T}}$ by
\begin{equation}
v_{0,\bar{T}} (\. ) = x_{0}b_{0}(\. ) + x_{\T } b_{n+1}(\. )
\end{equation}
The functions $b_{i}$, $0 \leq i \leq n+1$, can be chosen in $C^{1}$, which is only slightly
more restrictive than $H^{1} $ and more convenient for computations.
The space generated by the functions $\{ b_{i} \}_{i=0}^{n+1} \subset C^{1}(\bar{I}) $
will be denoted by $\EE$, so that $\EE \h \subsetneq \EE$.

Then, by letting the test function $\f $ run through the basis $\{ b_{i} \}_{i=1}^{n}$,
and asking that the weak form of the equation be satisfied by
$u_{ap}(\. ) = \sum u_{i}b_{i} (\. )$ and for each such $\f $, one can solve
the linear system with the coordinates $u_{i}$ as unknowns and with
$x_{\bar{T}}$ as parameter. This linear system reads
\begin{equation} \label{eq linear prob proj}
[\B ].\{ u \} = \{ \F \}
\end{equation}
Here, $[\B ] = [\B _{ij}]$ is the symmetric $n \times n $ matrix whose entries are
\begin{equation}
\B _{ij} = - \langle
\dot{b}_{i}(\. ) , \dot{b}_{j}(\. )
\rangle + k \langle
b_{i}(\. ) , b_{j}(\. )
\rangle , 1 \leq i,j \leq n
\end{equation}
the vector $\{ u \}$ is the $n \times 1$ vector of the unknown coefficients
\begin{equation}
u_{i} , 1 \leq i \leq n 
\end{equation}
and $\{ \F \} $ is the $n \times 1$ vector of the known forces,
\begin{equation}
\F _{i} = \langle
f(\. ) , b_{i}(\. )
\rangle -
 \langle
x_{0} \dot{b}_{0}(\. ) + x_{\T } \dot{b}_{n+1}(\. ) , \dot{b}_{i}(\. )
\rangle +
k \langle
x_{0}b_{0}(\. ) + x_{\T } b_{n+1}(\. ) , b_{i}(\. )
\rangle
\end{equation}
where $1 \leq i \leq n$. The final displacement, $x(\T ) = x_{\T }$ appears
here as a parameter and the solution $\{ u \}$ of the linear problem depends
on it.

The matrix $[ \B ] $ is invertible thanks to the choice of $\T $, which grants a
unique vector $\{ u \} = \{ u \} ( x_{\T } )$ satisfying eq. \eqref{eq linear prob proj}
Then the condition that $\dot{x}_{ap} (0) = \xx _{0}$, where
\begin{equation} \label{eqapprox sol}
x_{ap}(\. ) = \sum _{i=0}^{n+1} u_{i} b_{i}(\. ) \in \EE,
\end{equation}
is the approximate solution, allows to rid the linear system of eq. \eqref{eq linear prob proj}
of the dependence on $x_{\T }$.

The action of letting the test function run through a basis of the space
$\EH$ in which the approximate solution is sought does not allow
any arbitrary choices of parameters, and lets the geometry of the spaces
optimize the approximation. The only arbitrary choice consists in the choice
of the space $\mathcal{B}$ and in the construction of the basis $\{ b_{i} \}$
which should respectively be "as dense in $H^{1}_{0}$ as possible" and
facilitate calculations by satisfying good orthogonality conditions.

An important note is due at this point. In reality, the SDOF problem is presented as an
initial values problem. The weak form is better adapted to a boundary value problem.
Under the standing assumption on $\T$, initial and boundary value problems are shown to
be in a one-to-one correspondence for the \textit{exact solution}. When approximate
solutions are constructed, however, the dependence of the operator $\BB $ of eq.
\eqref{eq bound prob} on the excitation function $f$ comes into play. The correspondence
\begin{equation}
\dx _{0} \leftrightarrow x _{ap}(\T )= \BB _{x_{0},f_{ap},\bar{T}}(\dx _{0})
\end{equation}
depends, therefore, on the space $\EE$ in a
quite subtle way.

As is expected, $x _{ap}(\T ) \ra \BB _{x_{0},f,\bar{T}}(\dx _{0})$ as $x_{ap} \ra x$,
provided that $f_{ap} \ra f $ in $H^{-1}$. However, equality is in general not true
before passage to the limit and convergence needs to be shown,
since the solution $x_{ap}$ as constructed in eq. \eqref{eqapprox sol} satisfies the
initial conditions $x_{ap}(0) = x_{0}$ and $\dx_{ap}(0) = \dx _{0}$ and \textit{not} the
boundary condition $x_{ap}(\T) = x(\T)$.

For this reason, all statements concerning the construction of the
approximate solution will be stated in terms of initial values, and not in terms of
boundary conditions.

Concerning the choice of the space $\EE$, the
authors have tried to use "damped stationary waves", i.e. functions like
\begin{equation}
e^{-c t /2 }
\sin \left( \pi k \frac{t }{\bar{T}} \right) , k \in \N ^{*}
\end{equation}
for generating the space $\mathcal{B}$, coupled with cosines carrying the
boundary conditions. The attempt was unsuccessful, but a successful one
using $C ^{1} $ piecewise polynomial functions will be presented in part II of
the paper.


\subsection{The excitation function of the approximate solution}

Let $\EE$ be as in the previous paragraph and
$\{ b_{i} \}_{i=0}^{n+1} \subset C^{\infty }_{0} $ be a basis for $\EE$, let
$\EI \subset \EE$ be the set of functions in $\EE$ satisfying homogeneous
\textit{initial} conditions:
\begin{equation}
\EI = \{
b \in \EE, b(0) = \dot{b}(0) = 0
\}
\end{equation}
and suppose for convenience and without loss of generality that the first two functions are the
only two ones with non-zero derivative at $0$.\footnote{It is a standing assumption that
$\{ b_{i} \}_{i=1}^{n}$ satisfy homogeneous boundary conditions.}, so that
\begin{equation}
\begin{array}{r@{}l}
\EE _{ic} &= \R b_{0} \oplus \R b_{1} \\
\EE &= \EI \oplus \EE _{ic}
\end{array}
\end{equation}

Define also the linear mapping $\Phi $, associating to a displacement function
$x \in H^{1}$, the corresponding functional in $H^{-1}$ defining the lhs of the weak
formulation,
\begin{equation}
\Phi (x) = \Phi _{c} (x) = 
\left(
\f (\. ) \ra - \langle
\xx (\. ) , \dot{\f} (\. )
\rangle _{c}
+ k 
\langle
x (\. ) , \f (\. )
\rangle _{c}
\right)
\end{equation}
This mapping is continuous $H^{1} \ra H^{-1}$, and the following
proposition is obvious, given theorem \ref{thm weak form SDOF}.
\begin{prop} 
Under the assumptions of theorem \ref{thm weak form SDOF}, with homogeneous boundary
conditions, the mapping $\Phi$ a bijection $H^{1}_{0} \ra H^{-1}$. The function $x$ is
a solution to the problem if, and only if,
\[
\Phi _{c} (x) (u) = \langle e^{c\.}f, u \rangle _{H^{-1}} =
\langle F, \dot{u} \rangle_{L^{2}_{0}}, \forall u \in H^{1}_{0}
\]
\end{prop}

Define, also, the following spaces of excitation functions corresponding to the space $\EE$
and its relevant subspaces:
\begin{equation} \label{eqdef force spaces}
\begin{array}{r@{}l}
\FFF &= \FF(\EE)\\ 
\FFF _{00} &= \FF(\EI) \\ 
\FFF_{0}^{0} &= \FF( \EH ) \\ 
\FFF_{ic} &= \FF(  \EE _{ic})
\end{array}
\end{equation}
where $\FF$ is defined in eq. \eqref{eqdefF}.
All these spaces are subspaces of $H^{-1}$, that can be identified to subspaces of
$L^{2}_{0}$ by means of the $_{H^{-1}}\langle \. , \. \rangle _{H^{1}}$ duality via the
mapping $H^{-1} \ra L^{2}_{0}$
\begin{equation}
e^{c\.}f \mapsto F = \int e^{c\.}f
\end{equation}
as in lemma \ref{lemequiv f and F}. The shorthand in notation will be preserved when
subscripts are used, so that
\begin{equation}
F_{\#} = \int e^{c\.}f_{\#}
\end{equation}
whatever the subscript $\#$ might be.

Let, now, $x$ be the exact solution of the SDOF problem with rhs equal to $f$, and initial
conditions $x_{0}, \dx_{0}$. There exists
$\l = \l (x_{0}, \dx_{0})\in \R$, $|\l | \leq K \max (x_{0}, \dx_{0})$ where $K>0$ depends
only on the choice of $b_{0}$ and $b_{1}$, such that
\begin{equation}
x_{h} = x - x_{0}b_{0} - \l b_{1}
\end{equation}
satisfies homogeneous initial conditions. Call
\begin{equation}
f_{h} = \FF (x_{h}) = f - x_{0} \FF(b_{0}) - \l \FF ( b_{1})
\end{equation}
the corresponding excitation function. The same procedure can be applied to the approximate
solution $x_{ap}$ constructed in \S \ref{seccons app sol}, which satisfies the
strong SDOF equation with rhs $f_{ap}$, yielding $f_{ap,h}$.

With these notations the following proposition holds, whose proof follows directly from
the definitions of the respective functional spaces.
\begin{prop} \label{propexc func approx}
The excitation function $f_{ap}$ corresponding to the approximate solution $x_{ap}$
as in eq. \eqref{eqapprox sol} satisfies
\begin{equation}
\begin{array}{r@{}l}
f_{h} &= f - x_{0} \FF(b_{0}) - \l \FF ( b_{1})\\
f_{ap} &\in\FFF \\
f_{ap,h} &= f_{ap} - x_{0}   \FF(b_{0})
- \l  \FF ( b_{1}) \in \FFF_{00}\\
\pi _{\EE \h}e^{c\.} f_{ap, h} &= \pi _{\EE \h}e^{c\.} f_{h}
\end{array}
\end{equation}
where the projection $\pi _{\EH}$ is with respect to the $H^{-1}$ duality.
\end{prop}
The definition of the projection $\pi $ with respect to the duality
$_{H^{-1}}\langle \. , \. \rangle _{H^{1}}$ uses the injections
$H^{1} \hra L^{2} \hra H^{-1}$ where the space $L^{2} \equiv L^{2}_{0}$ is identified with
its dual. Thus, projection of $u \in H^{-1} $ along $\f \in H^{1}_{0 } $ with respect to
the $_{H^{-1}}\langle \. , \. \rangle _{H^{1}}$ duality amounts to projection of $U \in L^{2}$, an integral of $e^{c\.}u(\.)$, along $\dot{\f} \in L^{2}$:
\begin{equation}
\begin{array}{r@{}l}
\pi _{\R \f} e^{c\.} u &=
\frac{1}{\|\f\|_{H^{1}}^{2}}
\langle u , \f \rangle_{c}
\f
\\
&=
\frac{1}{\|\f\|_{H^{1}}^{2}}
\langle U , \dot{\f}\rangle_{0}
\f
\end{array}
\end{equation}

The proof follows from calculations similar to the ones proving theorem
\ref{thm weak form SDOF} and the construction of the approximate solution.

%

The following proposition is immediate.
\begin{prop}
The following estimate on $f_{ap,h}$ is true:
\begin{equation}
\|f_{ap,h} \|_{H^{-1}} \leq \|f \|_{H^{-1}} + K \max (x_{0},x_{\T})
\end{equation}
\end{prop}

From the properties of projections in Hilbert spaces follows directly the following
corollary to proposition \ref{propexc func approx}.

\begin{corollary}
The following estimate on the excitation function of the approximate solution holds
true:
\begin{equation}
\| e^{c\.}f_{ap, h} - \pi _{\EH}e^{c\.} f_{h}\|_{H^{-1}} \lesssim
\|\pi _{\EH} e^{c\.}f_{h}\|_{H^{-1}} \tan \theta _{nh}+
K \max (x_{0},x_{\T}) \sin \theta _{h}
\end{equation}
where $\theta _{nh}$ is the angle between the spaces $e^{c\.} \FFF _{00}$ and $e^{c\.}\EH$ in $H^{-1}$
and $\theta _{h}$ is the angle between the spaces $e^{c\.} \FFF _{ic}$ and $e^{c\.}\EH$ in $H^{-1}$.
\end{corollary}
Equivalently and more easily in terms of calculation, the angle $\theta _{nh}$ can be
obtained as the angle between the space formed by integrals of functions
$e^{c\. } u \in \FFF _{00}$ and derivatives of functions in $\EH$ in $L^{2}_{0}$. Similarly
for the angle $\theta _{h}$.

The subscripts in the angles $\theta$ stand for \textit{non-homogeneous} and
\textit{homogeneous} respectively, referring to the SDOF problem with $0$ initial
conditions and non-$0$ excitation function; and $0$ excitation function and non-$0$
initial conditions, respectively.

The choice of the space $\EE$ determines the geometry that is captured in the above
corollary: the better the corresponding spaces are aligned, i.e. the smaller the angles
$\theta$ are, the better the approximation is.

The following proposition and its corollary conclude the study of the convergence of the
numerical method in the general case.
\begin{prop} \label{propconv approx}
The error in the excitation force, $f_{er}$, given by
\begin{equation}
f_{er,h} = f_{h} - f_{ap,h} = (f_{h} - \pi_{\EH}f_{h}) + (\pi_{\EH}f_{h} - f_{ap,h})
\end{equation}
satisfies the following estimate
\begin{equation}
\|f_{er,h}\|_{H^{-1}} \leq \|f_{h} - \pi_{\EH}f_{h}\|_{H^{-1}} +
\|\pi _{\EH} f_{h}\|_{H^{-1}} \tan  \theta _{nh}+
K \max (|x_{0}|,|\dx_{0}|) \sin \theta _{h}
\end{equation}
\end{prop}
The proof follows immediately from the previous corollary and the triangle inequality.

\begin{corollary} \label{corconv approx}
The error in the excitation force for a given SDOF system can be bounded by
\begin{equation}
\|f_{er,h}\|_{H^{-1}} = O\left(\max
(
\|f_{h} - \pi_{\EH}f_{h}\|_{H^{-1}},
\|\pi_{\EH}f_{h}\|_{H^{-1}}\tan \theta _{nh},
\sin \theta _{h}
)
\right)
\end{equation}
\end{corollary}
In applications, the three terms will have to be estimated independently. The first
one stands for the error due to projection and it depends on how $f$ is placed with respect
to the space $\EH$. The second term depends on the relative position of the spaces
$\FFF_{0}^{0}$ and $\EH$, and the last one on the relative position of the spaces
$\R \Phi(e^{c\.} b_{0}) \oplus \R \Phi (e^{c\.} b_{1})$ and $\EH$.
For the approximation to converge, all three terms need to go to $0$ as
$\EE \ra \HH$. 

\subsection{Bounds for the approximation error} \label{secbounds theory}

In this paragraph, a method for obtaining bounds for the error of approximation is
presented under some assumptions on the $n$-dimensional space $\EH$, assuming
for convenience that the functions in the space $\EH$ are $C^{\infty } $ smooth.

%
%
%

To this end, call
\begin{equation}
\begin{array}{r@{}l}
\ddot{x} + c \dx + k x &= f \\
\ddot{x}_{ap} + c \dx _{ap} + k x _{ap} &= f_{ap} \\
\ddot{x}_{ap,h} + c \dx _{ap,h} + k x _{ap,h} &= f_{ap,h} \\
\ddot{x}_{er,h} + c \dx _{er,h} + k x  _{er,h}&= f _{er,h} = f_{h} - f _{ap,h} \\
\end{array}
\end{equation}
where $x_{ap,h} \in \EI$ satisfies homogeneous initial conditions, and $f_{ap,h}$
is obtained by substitution of $x_{ap,h}$ in the SDOF equation.
%
%
%
%

\subsubsection{Bounds for the error in displacement}

Firstly, the Duhamel representation, eq. \eqref{eq Duh}, and the calculations carried out
in eq. \eqref{eqsol Duh F} give the following proposition.
\begin{prop} \label{lemest duhamel}
Let $f \in H^{-1} $ and $F \in L^{2}$ an integral of $e^{c\. }f(\. )$. Then, for some
irrelevant constant $K>0$, the solution $x$ to the SDOF problem with homogeneous initial
conditions satisfies
\begin{equation}
\begin{array}{r@{}l}
\| x \|_{L^{\infty}} &\leq
K
\max \limits _{t \in I}e^{-ct/2}
\left(
\frac{1-e^{-ct }}{c}
\right)^{1/2}
\| F(\. )\|_{L^{2}_{0}(0,t)}
\leq
K
\| F(\. )\|_{L^{2}_{0}(I)}
\\
\| x \|_{L^{2}_{0}}&\leq
K\
e^{-c\T/2 }
\left(
1-e^{-c\T }
\right)^{1/2}
\| F(\. )\|_{L^{2}_{0}(I)}
\end{array}
\end{equation}
The constant $K$ depends on $\T $ and the parameters defining the SDOF system.
\end{prop}

\begin{proof}
The lemma amounts to estimating
\begin{equation}
x(t)
=
d(t)
 \int_{0}^{t} F(\tau  ) (d(\tau) -\frac{c}{2}s(\tau) ) \dd\tau
+s(t)
\int_{0}^{t} F(\tau  ) (\frac{c}{2}d(\tau) + \w _{d}^{2} s(\tau) ) \dd\tau
\end{equation}
in $L^{2}_{0}$ and in $L^{\infty }$. One can get directly that, for some irrelevant
constant $K>0$ depending on $c$ and $\w _{d}$,
\begin{equation}
\begin{array}{r@{}l}
|x(t)|
&\leq
Ke^{-ct/2}
\left(
(\int _{0}^{t} F(\tau )d(\tau) \dd\tau )^{2}
+
(\int _{0}^{t} F(\tau )s(\tau) \dd\tau )^{2}
\right)^{1/2}
\\
&\leq
Ke^{-ct/2}
\left(
\int _{0}^{t} d(\tau) ^{2}\dd\tau
+
\int _{0}^{t} s(\tau)^{2} \dd\tau 
\right)^{1/2}
(\int _{0}^{t} F(\tau )^{2}d(\tau) \dd\tau )^{1/2}
\\
&\leq
Ke^{-ct/2}
\left(
\frac{1-e^{-ct }}{c}
\right)^{1/2}
\| F(\. )\|_{L^{2}_{0}(0,t)}
\end{array}
\end{equation}
by repeated applications of H\"{o}lder's inequality. Direct estimation in $L^{2}$ gives
\begin{equation}
\begin{array}{r@{}l}
\|x(t)\|_{L^{2}_{0}(I)}^{2}
&\leq
K\|
\frac{e^{-ct }-e^{-2ct }}{c}
\|_{L^{1}_{0}(I)}
\| \| F(\. )\|_{L^{2}_{0}(0,t)} \|^{2}\|_{L^{\infty}(I,dt)}
\\
&=
Ke^{-c\T }(1-e^{-c\T })
\| F(\. )\|_{L^{2}_{0}(I)}^{2}
\end{array}
\end{equation}
which concludes the proof.
\end{proof}

Using lemma \ref{lemequiv f and F}, the following corollary is immediate.

\begin{corollary} \label{corconv x in l2}
The error in displacement converges to $0 $ in $L^{2}_{0}$ at the same rate as
$\|f _{er,h}\| _{H^{-1}}$:
\begin{equation}
\| x _{er,h}\|_{L^{2}_{0}} \leq
C \|f _{er,h} \| _{H^{-1}}
\end{equation}
\end{corollary}

Under the standing assumptions on $\EE$, bounds for the errors of the approximate solution
of the SDOF problem can be obtained by applying the following proposition.

\begin{prop} \label{propbounds disp}
Let $x $ be the solution to a weak SDOF problem, $f_{h}$ the excitation function
corresponding to homogeneous initial conditions, see proposition \ref{propexc func approx},
and $F _{h}  \in L^{2}_{0}$ an integral of $e^{c\.}f_{h}(\. )$. Then, the error of
approximation $x_{er}$ satisfies the following estimates:
\begin{equation}
\begin{array}{r@{}l}
\|F _{h}  \|_{L^{2}_{0}} &\leq \|f \|_{H^{-1}} + K\max (x_{0},\dx _{0})
\\
\| x_{er,h} \|_{L^{\infty}}
&\leq
K
\max \limits _{t \in I}
e^{-ct/2}
\left(
\frac{1-e^{-ct }}{c}
\right)^{1/2}
\| F_{h,er}(\. )\|_{L^{2}_{0}(0,t)}
\left(
\| d_{t}(\. )\|_{L^{2}_{0}(0,t)}^{2}
+
\| s_{t}(\. )\|_{L^{2}_{0}(0,t)}^{2}
\right)^{1/2}
\\
&\leq
K
\| F_{h,er}(\. )\|_{L^{2}_{0}(I)}
\left(
\| d_{h}(\. )\|_{L^{2}_{0}(I)}^{2}
+
\| s_{h}(\. )\|_{L^{2}_{0}(I)}^{2}
\right)^{1/2}
\\
&\leq
K
\| F_{h,er}(\. )\|_{L^{2}_{0}(I)}
\\
\| x_{er,h} \|_{L^{2}_{0}} &\leq
K
e^{-c\T/2}
\left(
1-e^{-c\T }
\right)^{1/2}
\| F_{h,er}(\. )\|_{L^{2}_{0}(I)}
\\
| x_{er,h}(\T )| &\leq
K
e^{-c\T/2} \|F_{er,h}(\.  )\|_{L^{2}_{0}}
\left(
 \|d_{er,h}(\.) \|^{2}_{L^{2}_{0}} +
\|s_{er,h}(\.) \|^{2}_{L^{2}_{0}}
\right)^{1/2}
\end{array}
\end{equation}
\end{prop}
In the statement of the proposition,
\begin{equation}
\begin{array}{r@{}l}
F_{h,er} (\. ) &= \pi _{(\dot{\EE} _{0}^{0})^{\perp}}  F_{h}(\. )\\
&=
F_{h}(\. ) - \pi _{\dot{\EE} _h}  F_{h}(\. )
\\
d_{er,h} (\. ) &= \pi _{(\dot{\EE}_{0}^{0})^{\perp}}  d(\. )\\
&=
d(\. ) - \pi _{\dot{\EE}_{0}^{0}}  d(\. )
\\
s_{er,h} (\. ) &= \pi _{(\dot{\EE}_{0}^{0})^{\perp}}  s(\. )\\
&=
s(\. ) - \pi _{\dot{\EE}_{0}^{0}}  s(\. )
\end{array}
\end{equation}

The errors of approximation to the fundamental solutions $d_{er}(\. )$ and $s_{er} (\. )$
come from projections on the space $\dot{\EE} _{0}^{0}$ of (the derivatives of) test functions.
These errors naturally depend on the choice of the space $\EE $, and therefore on the
specifics of each implementation.
\begin{proof}
The proof uses the same calculations as the one of proposition \ref{lemest duhamel}, along
with proposition \ref{propexc func approx}, and the relation between
$H^{1}_{0}$-$H^{-1}$ duality and the $L^{2}$ scalar product of the proof of lemma
\ref{lemequiv f and F}. The last inequality follows by substitution $t = \T$:
\begin{equation}
\begin{array}{r@{}l}
x(\T) &= d(\T)
 \int_{0}^{\T} F(\tau  ) (d(\tau) -\frac{c}{2}s(\tau) ) \dd\tau +
s(\T)
\int_{0}^{\T} F(\tau  ) (\frac{c}{2}d(\tau) + \w _{d}^{2} s(\tau) ) \dd\tau
\\
&= d(\T)
 \langle F(\.  ), d(\.) -\frac{c}{2}s(\.) \rangle_0 +
s(\T)
\langle F(\.  ), \frac{c}{2}d(\.) + \w _{d}^{2} s(\.) \rangle_0
\end{array}
\end{equation}
Using the fact that $|d(\.)|^{2}+ \w_d^{2}|s(\.)|^{2} = e^{-c\.}$, one gets
directly that
\begin{equation}
\begin{array}{r@{}l}
| x(\T) |
&\leq
e^{-c\T/2} \left(
\langle F(\.  ), d(\.) -\frac{c}{2}s(\.) \rangle_0^{2} +
\langle F(\.  ), \frac{c}{2\w _{d}}d(\.) + \w _{d} s(\.) \rangle_0^{2}
\right)^{1/2}
\\
&\leq
e^{-c\T/2} \|F(\.  )\|_{L^{2}_{0}}
\left(
 \|d(\.) -\frac{c}{2}s(\.) \|^{2}_{L^{2}_{0}} +
\|\frac{c}{2\w _{d}}d(\.) + \w _{d} s(\.) \|^{2}_{L^{2}_{0}}
\right)^{1/2}
\\
&\leq
e^{-c\T/2} \|F(\.  )\|_{L^{2}_{0}}
\left(
 \|d(\.) -\frac{c}{2}s(\.) \|^{2}_{L^{2}_{0}} +
\|\frac{c}{2\w _{d}}d(\.) + \w _{d} s(\.) \|^{2}_{L^{2}_{0}}
\right)^{1/2}
\end{array}
\end{equation}
which implies the desired result as soon as the calculation is carried out for
$F_{er,h}$, and the projections of the fundamental solutions are substituted for the
exact fundamental solutions.
\end{proof}
It can be noted that further calculations can provide the bound
\begin{equation}
| x(\T) |
\leq
e^{-c\T/2} (1 + \frac{c^{2}}{4\w_{d}^{2}}\frac{1-e^{-c\T}}{c} )\|F(\.  )\|_{L^{2}_{0}}
\end{equation}
for the exact solution, but the expression as provided in the statement of the
proposition will be more useful in practice, since it comprises the projections of
the fundamental solutions, resulting in improved bounds.

This proposition provides an easier way of estimating the errors, since it uses function
spaces instead of distribution ones, and measures the error of approximation by the error
of approximation to the two fundamental solutions of the problem, which do not depend on
the excitation function.


\subsubsection{Bounds for the error in velocity}

The formulas of eq. \eqref{eqder fund sol} imply directly that the terms of the first line
are of the order of $x(t)$, while that of the second line should be expected to be of a
greater order, as $F(\. )$ is a priori only $L^{2}$ regular, while its integral is in
$H^{1}$. The following proposition and its corollaries, the counterparts of
proposition \ref{lemest duhamel} and corollary \ref{corconv x in l2} for the velocity,
establish this fact.
\begin{prop} \label{lemest duhamel dx}
Let $f \in H^{-1} $ and $F \in L^{2}$ an integral of $e^{c\. }f(\. )$. Then, for some
irrelevant constant $K>0$,
\begin{equation}
\| \dx \|_{L^{2}_{0}} \leq
K
\left(
\| x \|_{L^{2}_{0}}
+
\| e^{-c\. } F(\. )\|_{L^{2}_{0}(I)}
\right)
\end{equation}
The constant $K$ depends on $\T $ and the parameters defining the SDOF system.
%
\end{prop}
The following corollary follows directly from the proposition.
\begin{corollary} \label{corconv dx in l2}
The error in displacement converges to $0 $ in $L^{2}_{c}$ at the same rate as
$\|f _{h,er}\| _{H^{-1}}$:
\begin{equation}
\| \dx _{er}\|_{L^{2}_{c}} \leq
K \|f _{h,er} \| _{H^{-1}}
\end{equation}
\end{corollary}

The following proposition is the counterpart to prop. \ref{propbounds disp}, bounding the
error in the approximation of the velocity in a similar manner. However, there is no
$L^{\infty}$ estimation, since a priori $F$ and therefore $\dx $ are only $L^{2}$ and not
$L^{\infty}$.

For the estimation of the velocity at time $t=\T$, the following calculation will be
necessary. Taking the derivative of eq. \eqref{eqsol Duh F} gives
\begin{equation}
\begin{array}{r@{}l}
\dx(t)
&=
\dot{d}(t)
 \int_{0}^{t} F(\tau  ) (d(\tau) -\frac{c}{2}s(\tau) ) \dd\tau
+\dot{s}(t)
\int_{0}^{t} F(\tau  ) (\frac{c}{2}d(\tau) + \w _{d}^{2} s(\tau) ) \dd\tau
\\
&\phantom{d(t) \int_{0}^{t} F(\tau  ) (d(\tau)}+
(d(t)^{2}+\w _{d}^{2}s(t)^{2})F(t)
\\
&=
\dot{d}(t)
 \int_{0}^{t} F(\tau  ) (d(\tau) -\frac{c}{2}s(\tau) ) \dd\tau
+\dot{s}(t)
\int_{0}^{t} F(\tau  ) (\frac{c}{2}d(\tau) + \w _{d}^{2} s(\tau) ) \dd\tau
\\
&\phantom{d(t) \int_{0}^{t} F(\tau  ) (d(\tau)}
+ e^{-ct}F(t)
\end{array}
\end{equation}

\begin{prop} \label{propbounds speed}
Let $x $ be the solution to a weak SDOF problem, $f_{h}$ the excitation function
corresponding to homogeneous initial conditions, see proposition \ref{propexc func approx},
and $F _{h}  \in L^{2}_{0}$ an integral of $e^{c\.}f_{h}(\. )$. Then, the error of
approximation $\dx_{er}$ satisfies the following estimates:
\begin{equation}
\begin{array}{r@{}l}
\| \dx_{er,h} \|_{L^{2}_{0}} &\leq
K
\left(
\frac{1-e^{-2c\T }}{2c}
\right)^{1/2}
\| F_{er,h}(\. )\|_{L^{2}_{0}(I)}
+ K \| x_{er,h} \|_{L^{2}_{0}}
\\
| \dx_{er,h}(\T )| &\leq
e^{-c\T /2}
\| F_{er,h}(\. )\|_{L^{1}(I)}
+ K |x_{er,h}(\T )|
\\
&\leq
e^{-c\T /2} \sqrt{\T}
\| F_{er,h}(\. )\|_{L^{2}(I)}
+ K |x_{er,h}(\T )|
\end{array}
\end{equation}
\end{prop}

The proof is similar to that of proposition \ref{propbounds disp}.
To recall, $f$ is supposed $L^{2}$ close to the boundary, which makes $F(\T)$
unambiguously defined. The quantity $\| F_{h,er}(\. )\|_{L^{1}(I)}$ is also well
defined thanks to the good boundary behavior of $f$.

\begin{proof}
Concerning the estimate at time $t=\T$, one gets
\begin{equation}
\begin{array}{r@{}l}
\dx(\T) &= e^{-ch}F(\T) +
(d(\T) -\frac{c}{2}s(\T) )
\int_{0}^{h} F(\tau  ) (\frac{c}{2}d(\tau) + \w _{d}^{2} s(\tau) ) \dd\tau
\\
&\phantom{\frac{c}{2}s(t) )
\int_{0}^{t} F(\tau  ) (\frac{c}{2}d(\tau)}
-(\frac{c}{2}d(\T) + \w _{d}^{2} s(\T)
 \int_{0}^{h} F(\tau  ) (d(\tau) -\frac{c}{2}s(\tau) ) \dd\tau
\\
&= e^{-ch}F(\T) +
(d(\T) -\frac{c}{2}s(\T) )
\langle  F(\.  ) ,(\frac{c}{2}d(\.) + \w _{d}^{2} s(\.) ) \rangle_0
\\
&\phantom{(\frac{c}{2}s(t) )
 F(\.  ) (\frac{c}{2}d(\.)}
-(\frac{c}{2}d(\T) + \w _{d}^{2} s(\T))
 \langle  F(\.  ), (d(\.) -\frac{c}{2}s(\.) ) \rangle_0
\end{array}
\end{equation}
so that
\begin{equation}
\begin{array}{r@{}l}
|\dx(\T)| &\leq
|e^{-ch}F(\T)| +
\left(
(d(\T) -\frac{c}{2}s(\T) )^{2}
\|  F(\.  )\|^{2}_{L^{2}_{0}} \| (\frac{c}{2}d(\.) + \w _{d}^{2} s(\.) )\|^{2}_{L^{2}_{0}}
\right.
\\
&\phantom{( F(\.  ) (\frac{c}{2}d(\.)}
\left. 
+(\frac{c}{2}d(\T) + \w _{d}^{2} s(\T))^{2}
 \|  F(\.  )^{2}_{L^{2}_{0}} \| (d(\.) -\frac{c}{2}s(\.) ) ^{2}_{L^{2}_{0}}
\right)^{1/2}
\end{array}
\end{equation}
from which the estimate of the proposition follows immediately.
\end{proof}

\subsubsection{Some energy estimates}

The conservation (cf. \S \ref{seccons law}) and dissipation (cf. \S \ref{seccdis law})
estimates, together with the above proposition and proposition \ref{lemest duhamel} provide the
following corollaries.
\begin{corollary} \label{corconv H1}
Let $x $ and $x_{er}$ as in prop. \ref{propbounds disp}. Then,
\begin{equation}
\| \dx _{er} \|_{L^{2}_{0}} + \| x _{er} \|_{L^{2}_{0}} \leq
K \|F_{h,er}\|_{L^{2}_{0}} 
\end{equation}
where the constant $C$ depends on the space $\EE $, the parameters of the SDOF system
and $\T $.
\end{corollary}


Concerning the convergence of the velocity, $\dx$, lemma \ref{corconv H1} already provides bounds
in $\HH$. In order to obtain uniform convergence, $H^{2}$ bounds are needed. They are
provided by the following obvious corollary.
\begin{corollary} \label{corconv H2}
Let $f _{h,er}$ and $x_{er}$ be as in prop. \ref{propbounds disp}. Then,
\begin{equation}
\| \ddot{x} _{er} \|_{H^{-1}} \leq
K \|F_{h,er}\|_{L^{2}_{0}}
\end{equation}
\end{corollary}
\begin{proof}
The acceleration satisfies the equation
\begin{equation}
\ddot{x} = f - c \dx - kx
\end{equation}
Estimation by duality, the triangle inequality and corollary \ref{corconv H1} conclude the
proof.
\end{proof}

\subsection{Convergence}

Summing up, the following theorem holds true.
\begin{theorem} \label{thmconv}
Suppose a sequence of spaces $\EE _{n} $ is given and such that
$\cup \dot{\EE} _{n} = L^{2} \mod \R$.
Then, the sequence of corresponding approximate solution $x^{(n)}_{ap} \in \EE _{n}$
converges to $x$ in $H^{1}$:
\begin{equation}
\|x^{(n)}_{er}\|_{H^{1}} = O(\|F_{h,er}^{(n)}\|_{L^{2}_{0}}) \ra 0
\end{equation}
In particular, since $\T <\infty$,
\begin{equation}
\|x^{(n)}_{er}\|_{L^{\infty}} = O(\|F_{h,er}^{(n)}\|_{L^{2}_{0}}) \ra 0
\end{equation}
\end{theorem}
\begin{proof}
The theorem follows from corollaries \ref{corconv x in l2} and \ref{corconv dx in l2}.
The uniform convergence for the displacement follows from Sobolev the injection
$\HH_{0} \hra L^{\infty}$ which holds for finite intervals.
\end{proof}

Under additional assumptions on the excitation function, and by invoking the energy
estimates, the following stronger theorem can be proved.

\begin{theorem} \label{thmconv smooth}
Suppose a sequence of spaces $\EE _{n} $ is given and such that $\cup \EE _{n} =
 L^{2} \mod \R$.
Suppose, moreover, that $F_{er}^{(n)}$ stays bounded in the $\HH$ norm. Then, the
sequence of corresponding approximate solution $x^{(n)}_{ap} \in \EE _{n}$
converges to $x$ in $H^{2}$:
\begin{equation}
\|x^{(n)}_{er}\|_{H^{2}} = O(\|F_{h,er}^{(n)}\|_{L^{2}_{c}}) \ra 0
\end{equation}
In particular, since $\T <\infty$,
\begin{equation}
\|x^{(n)}_{er}\|_{L^{\infty}} +  \|\dx ^{(n)}_{er}\|_{L^{\infty}}
= O(\|F_{h,er}^{(n)}\|_{L^{2}_{c}}) \ra 0
\end{equation}
\end{theorem}
\begin{proof}
The theorem follows from propositions \ref{propbounds disp} and \ref{propbounds speed}.
\end{proof}
Summing up, suppose that one can consider a sequence of spaces $\EE _{n} $
satisfying the assumptions of this section and such that $\cup \EE _{n} = L^{2}_{c} $.
Then, the sequence of corresponding approximate solution $x^{(n)}_{ap} \in \EE _{n}$
converges to $x$ in $H^{2}$ as long as $f -f_{ap}$ stays bounded in the $\HH$ norm. This convergence
implies that
\begin{equation}
\begin{array}{r@{}l}
x^{(n)}_{ap} &\ra x \text{ uniformly} \\
\dx^{(n)}_{ap} &\ra \dx \text{ uniformly} \\
\ddot{x}^{(n)}_{ap} &\ra \ddot{x} \text{ in } L^{2}_{c}
\end{array}
\end{equation}
and the rates of convergence are given by corollary \ref{corconv H2}.


Similarly one can compute
\begin{equation}
\begin{array}{r@{}l}
\dx _{er}(t) &= -\frac{c}{{2\omega _d }}\int_0^t f _{er}(\tau )e^{ - c (t
- \tau )/2} \sin [\omega _d (t - \tau )]\dd\tau
+\\
&\phantom{\int_0^t f _{er}(\tau )e^{ - c (t
- \tau )/2} \sin [\omega _d }
+\int_0^t f _{er}(\tau )e^{ - c (t
- \tau )/2} \cos [\omega _d (t - \tau )]\dd\tau
\end{array}
\end{equation}
which implies that
\begin{equation}
\begin{array}{r@{}l}
|\dx _{er}(\T )| & \leq -\frac{ce^{ - c \T /2}}{{2\omega _d }}\int_{0}^{\T } f _{er}(\tau )
e^{ c \tau /2} \sin [\omega _d (\T - \tau )]\dd\tau
+\\
&\phantom{\int_{0}^{\T } f _{er}(\tau )e^{  c \tau /2} \sin [\omega _d }
+e^{ - c \T /2}\int_{0}^{\T } f _{er}(\tau )e^{ c \tau /2}
\cos [\omega _d (t - \tau )]\dd\tau
\end{array}
\end{equation}


\section{Conclusion}

This is the first of three parts of a work on the weak formulation of the SDOF and MDOF
problems. The weak formulation of both problems has been obtained, and a framework
for obtaining numerical methods using the weak formulation has been laid out. The
particular instances of the numerical methods using the weak formulation depend on the
choices of bases for the functional space $\HH $. The framework presented in this
work includes results about the convergence properties of the particular instances.

In part two of the work, such an instance will be presented, using Bernstein polynomials
as basis for $\HH $, and the convergence results will be specified in this context,
proving fast rates of convergence that are outside the scope of traditional step-wise
methods.

\bibliography{aomsample}

\providecommand{\bysame}{\leavevmode\hbox to3em{\hrulefill}\thinspace}
\providecommand{\noopsort}[1]{}
\providecommand{\mr}[1]{\href{http://www.ams.org/mathscinet-getitem?mr=#1}{MR~#1}}
\providecommand{\zbl}[1]{\href{http://www.zentralblatt-math.org/zmath/en/search/?q=an:#1}{Zbl~#1}}
\providecommand{\jfm}[1]{\href{http://www.emis.de/cgi-bin/JFM-item?#1}{JFM~#1}}
\providecommand{\arxiv}[1]{\href{http://www.arxiv.org/abs/#1}{arXiv~#1}}
\providecommand{\doi}[1]{\url{http://dx.doi.org/#1}}
\providecommand{\MR}{\relax\ifhmode\unskip\space\fi MR }
\providecommand{\MRhref}[2]{%
  \href{http://www.ams.org/mathscinet-getitem?mr=#1}{#2}
}
\providecommand{\href}[2]{#2}
\begin{thebibliography}{Rud87}

\bibitem[AM78]{AbrahamMechanics}
\bgroup\scshape{}R.~Abraham\egroup{} and \bgroup\scshape{}J.~Marsden\egroup{},
  \emph{Foundations of Mechanics}, \emph{AMS Chelsea publishing}, AMS Chelsea
  Pub./American Mathematical Society, 1978. Available at
  \url{https://books.google.com.br/books?id=4Y-ownk6ilsC}.

\bibitem[dB78]{dBSpl}
\bgroup\scshape{}C.~de~Boor\egroup{}, \emph{A practical guide to splines},
  \emph{Applied Mathematical Sciences} \textbf{27}, Springer-Verlag, New
  York-Berlin, 1978. \mr{507062}.

\bibitem[Br{\'e}83]{BrezisAnFonc}
\bgroup\scshape{}H.~Br{\'e}zis\egroup{}, \emph{Analyse fonctionnelle:
  th{\'e}orie et applications}, \emph{Collection Math{\'e}matiques
  appliqu{\'e}es pour la ma{\^\i}trise}, Masson, 1983. Available at
  \url{https://books.google.com.br/books?id=icO4ngEACAAJ}.

\bibitem[CQ82]{ApproximationLegSob}
\bgroup\scshape{}C.~Canuto\egroup{} and
  \bgroup\scshape{}A.~Quarteroni\egroup{}, Approximation results for orthogonal
  polynomials in sobolev spaces,  \emph{Mathematics of Computation} \textbf{38}
  (1982), 67--86.

\bibitem[CP75]{clough1975dynamics}
\bgroup\scshape{}R.~Clough\egroup{} and \bgroup\scshape{}J.~Penzien\egroup{},
  \emph{Dynamics of structures}, \emph{Dynamics of Structures Ray W. Clough,
  Joseph Penzien}, McGraw-Hill, 1975. Available at
  \url{https://books.google.com.br/books?id=UdxRAAAAMAAJ}.

\bibitem[Dah56]{Dahl56}
\bgroup\scshape{}G.~Dahlquist\egroup{}, Convergence and stability in the
  numerical integration of ordinary differential equations,  \emph{MATHEMATICA
  SCANDINAVICA} \textbf{4} (1956), 33--53. \doi{10.7146/math.scand.a-10454}.
  Available at \url{http://www.mscand.dk/article/view/10454}.

\bibitem[Dah63]{Dahl63}
\bgroup\scshape{}G.~G. Dahlquist\egroup{}, A special stability problem for
  linear multistep methods,  \emph{BIT Numerical Mathematics} \textbf{3}
  (1963), 27--43. \doi{10.1007/BF01963532}.  Available at
  \url{http://dx.doi.org/10.1007/BF01963532}.

\bibitem[Far00]{FAROUKI2000145}
\bgroup\scshape{}R.~T. Farouki\egroup{}, Legendre–bernstein basis
  transformations,  \emph{Journal of Computational and Applied Mathematics}
  \textbf{119} (2000), 145--160.
  \doi{https://doi.org/10.1016/S0377-0427(00)00376-9}.  Available at
  \url{https://www.sciencedirect.com/science/article/pii/S0377042700003769}.

\bibitem[HS74]{HirschSmaleODE}
\bgroup\scshape{}M.~W. Hirsch\egroup{} and \bgroup\scshape{}S.~Smale\egroup{},
  \emph{Differential equations, dynamical systems, and linear algebra},
  Academic Press New York, 1974 (English). Available at
  \url{http://www.loc.gov/catdir/toc/els031/73018951.html}.

\bibitem[Hug87]{Hugh87}
\bgroup\scshape{}T.~J.~R. Hughes\egroup{}, \emph{The finite element method :
  linear static and dynamic finite element analysis}, Englewood Cliffs, N.J.
  Prentice-Hall International, 1987. Available at
  \url{http://opac.inria.fr/record=b1086028}.

\bibitem[Kac38]{Kac1938}
\bgroup\scshape{}M.~Kac\egroup{}, Une remarque sur les polynomes de ms
  bernstein,  \emph{Studia Mathematica} \textbf{7} (1938), 49--51.

\bibitem[Lag88]{Lager88}
\bgroup\scshape{}P.~A. Lagerstrom\egroup{}, \emph{Matched asymptotic
  expansions}, \emph{Applied Mathematical Sciences} \textbf{76},
  Springer-Verlag, New York, 1988, Ideas and techniques. \mr{958913}.
  \doi{10.1007/978-1-4757-1990-1}.  Available at
  \url{http://dx.doi.org/10.1007/978-1-4757-1990-1}.

\bibitem[Rud76]{RudMathAn}
\bgroup\scshape{}W.~Rudin\egroup{}, \emph{Principles of mathematical analysis},
  3d ed. ed., McGraw-Hill New York, 1976 (English). Available at
  \url{http://www.loc.gov/catdir/toc/mh031/75017903.html}.

\bibitem[Rud87]{RudR&C}
\bgroup\scshape{}W.~Rudin\egroup{}, \emph{Real and complex analysis}, third
  ed., McGraw-Hill Book Co., New York, 1987. \mr{924157 (88k:00002)}.

\end{thebibliography}
\bibliographystyle{aomalpha}

\end{document}